\documentclass[twoside,12pt]{article}
\usepackage{amsmath}
\usepackage{amssymb}
\usepackage{amscd}
\usepackage{amsthm}%+

\numberwithin{equation}{section}

\theoremstyle{plain}
\newtheorem{theorem}{Theorem}[section]

\newtheorem{Le}[theorem]{Lemma}

\theoremstyle{remark}

\theoremstyle{definition}

\newcommand{\Z}{\mathbb Z}
\newcommand{\R}{\mathbb R}

\newcommand{\h}{\mathbb H}

\def\la{\lambda}

\newcommand{\vol}{\operatorname{Vol}}

\def\ga{\gamma}

\def\ra{\rightarrow}

\def\e{\emph}
\def\i{\infty}
\def\p{\partial}
\def\b{\begin}
\newcommand{\ol}{\overline}

\begin{document}

\title{%\flushleft
{Large scale geometry of   negatively curved   $\R^n \rtimes \R$}}
\author{Xiangdong Xie}
\date{  }

\maketitle

\begin{abstract}
We classify all negatively curved $\R^n \rtimes \R$ up to
quasiisometry.    We  show that all
 quasiisometries between  such manifolds
  (except
 when  they   are biLipschitz to the  real
hyperbolic spaces)
  are almost similarities.  
 We prove these results by studying the quasisymmetric
 maps on the ideal boundary of these manifolds.

\end{abstract}

%\noindent
%{\bf{Mathematics Subject Classification(2000).}} Primary 53C23, 51F99, 57M20;
% Secondary  53C70, 57M60.
% \newline

%\keywords    test  \endkeywords
%\subjclass  51  \endsubjclass

%\tableofcontents

%\vspace{3mm} \noindent
 {\small {\bf{Mathematics Subject
Classification (2000).}} 20F65,  30C65, 53C20.
%20F65. %57M20, 20F67, 20E07.}

%30C65  (1991-now) Quasiconformal mappings in ${\bf R}^n$, other generalizations
%20F34 Fundamental groups and their automorphisms
%20F65 Geometric group theory
%20F67 Hyperbolic groups and nonpositively curved groups
%20F99 None of the above, but in this section
%20E07 Subgroup theorems; subgroup growth
%53C20 Global Riemannian geometry, including pinching
%53C23 Global topological methods
%54F45 Dimension theory
%57M20 Two-dimensional complexes
%57M60 Group actions in low dimensions
%(1) (2) (a) (b) (i) (ii)
%20F69 Asymptotic properties of groups

%$\{x_i\}_{i=1}^\i$  converges to $\xi\in \ol{X}$
% an index two subgroup

               %\vspace{3mm} \noindent {\small {\bf{Key words.}} }

%hyperbolic element, parabolic element,
%quasi-convex.}

\tableofcontents

\setcounter{section}{0} \setcounter{subsection}{0}

\section{Introduction}\label{s0}

In this paper we  study quasiisometries between negatively curved
    solvable Lie groups of the form $\R^n\rtimes \R$  and quasisymmetric maps between  their ideal boundaries.

Given an $n\times n$  matrix $A$, we let $G_A$ be the semi-direct
product $\R^n\rtimes_A\R$, where $\R$ acts on $\R^n$ by
$(t,x)\mapsto e^{tA}x$.   Then $G_A$ is a solvable Lie
             group.

Let $G_A$ be equipped with any left invariant Riemannian metric
  such that the $\R$ direction is perpendicular to the $\R^n$
  factor.
  %determined by the standard inner product at the origin  of
   %$G_A=\R^n\times \R$.
When $A=I_n$,  $G_A$ is isometric to
$\h^{n+1}$. More generally,
             if  the eigenvalues of $A$ all have positive real parts,
             then  it follows from Heintze's results \cite{H} that $G_A$  is Gromov hyperbolic.
                %has  negative sectional curvature.
Hence $G_A$ has a well defined ideal boundary $\p G_A$.
  The ideal boundary $\p G_A$ can be naturally identified with (the
  one-point compactification of) $\R^n$. On the ideal boundary
   $\R^n$, there is a parabolic visual (quasi)metric  $D_A$, which
   is invariant under Euclidean translations and admits a family of
   dilations $\{\lambda_t=e^{tA}\}$.
        See Section \ref{metriconb} for more
     details.

Given an $n\times n$ matrix  $A$, the \e{real part Jordan form}  of
  $A$  is obtained from the Jordan form of $A$ by replacing each
  diagonal entry with  its real part and reordering  to make it
  canonical.  Notice that the real part Jordan form is different from
  the real Jordan form and the absolute Jordan form.  It is related
  to the  absolute Jordan form through matrix exponential.

     Here are  the main  results of the paper.
 See Theorem \ref{main2text} for a more precise statement of Theorem
 \ref{main2}.  Also see Section \ref{prelimi} for basic definitions.

     \b{Th}\label{main1}
Let $A$ and $B$ be  $n\times n$ matrices whose eigenvalues all have
positive real parts.  Then $(\R^n,  D_A)$ and $(\R^n, D_B)$ are
quasisymmetric if and only if there is some $s>0$ such that $A$ and
$sB$ have  the same real part  Jordan form.

     \end{Th}

     \b{Th}\label{main2}
Let $A$   and $B$ be   $n\times n$ matrices whose eigenvalues all
have positive real parts.
  Denote by  $\la_1$ and $\mu_1$  the smallest real parts of the eigenvalues
  of $A$ and $B$ respectively,  and   set $\epsilon=\la_1/\mu_1$.  %Let $k_A$ be the number of distinct
  %real parts of the eigenvalues of $A$.
  If  the real part Jordan form of $A$ is not a multiple of the identity
matrix $I_n$,  then  for every
   %$\eta$-
  quasisymmetric map $F:(\R^n,  D_A)  \ra (\R^n, D_B)$,
     %there is some $\epsilon>0$ such that
   the map
        $F:(\R^n,  D^\epsilon_A)  \ra (\R^n, D_B)$
    is  biLipschitz.
     % a  $K$-quasisimilarity, where $K$ depends only on  $\eta$,
     %$A$, $B$ and the metrics on $G_A$, $G_B$.
  %       $\epsilon$ and $k_A$.

     \end{Th}

When $A=I_n$, the manifold $G_A$ is isometric to  the real
hyperbolic space $\h^{n+1}$.  In this case, the ideal boundary is
$\R^n$ with the Euclidean metric, and hence the claim in
   Theorem   \ref{main2}   fails: there are non-biLipschitz quasiconformal maps in
     the Euclidean space
   $\R^n$.  More generally, if the real part Jordan form of  $A$
  is a multiple of $I_n$,  then it follows from the result of
  Farb-Mosher (see Section 2) that $(\R^n, D_A)$ is biLipschitz to
   $(\R^n, |\cdot|^\epsilon)$, where $|\cdot|$ denotes the Euclidean
   metric and $\epsilon>0$ is some constant.  Hence the claim in
   Theorem  \ref{main2}  also fails.

There are several consequences of the main results.

  Recall that  two   geodesic Gromov hyperbolic spaces admitting
  cocompact isometric group actions
    %Hadamard manifolds with pinched negative curvature
   are quasiisometric if and only if
  their
 ideal boundaries are quasisymmetric with respect to the visual
 metrics, see \cite{Pa} or \cite{BS}.  Hence
   Theorem \ref{main1}  yields the quasiisometric
 classification of all  negatively   curved  $\R^n\rtimes \R$.

\b{Cor}\label{c1} {Let $A$ and $B$ be  $n\times n$ matrices whose
eigenvalues all have positive real parts.  Then $G_A$ and $G_B$ are
quasiisometric if and only if there is some $s>0$ such that $A$ and
$sB$ have the same real part Jordan form.}

\end{Cor}

The next three results are consequences of  Theorem \ref{main2}.
%\cite{SX} for  detailed proofs.

A map $f: X\ra Y$ between two metric spaces is called an \e{almost
similarity}  if  there
 are
 constants $L>0$ and $C\ge 0$ such that
   $L\, d(x_1,  x_2)-C\le d(f(x_1), f(x_2))\le L\, d(x_1,  x_2)+C$ for all
    $x_1,  x_2\in  X$   and $d(y, f(X))\le C$ for all $y\in Y$.

\b{Cor}\label{c2} {Let $A$   and $B$ be   $n\times n$ matrices whose
eigenvalues all have positive real parts.  Suppose the real part
Jordan form of $A$ is not a multiple of the identity matrix $I_n$.
 Then every quasiisometry
  $f: G_A\ra G_B$
 is an almost similarity.}
  %That is, there
 %are
 %constants $L>0$ and $C\ge 0$ such that
  % $L\, d(x,y)-C\le d(f(x), f(y))\le L\, d(x,y)+C$ for all
   % $x,y\in G_A$.}

\end{Cor}

    We view  the canonical projection  $h_A: G_A=\R^n\times \R\ra \R$
    as the height function for $G_A$.
      Let
$A$ and $B$ be two $n\times n$ matrices.  %  whose eigenvalues have
   %positive real parts.
 A  quasiisometry  $f:
G_A\ra G_B$ is \e{height-respecting} if it maps the fibers of
  $h_A$   to
within uniformly bounded Hausdorff distance from the fibers of
  $h_B$.

%there is a map $h: \R\ra \R$ and a constant $C\ge 0$ such that the
%Hausdorff distance between $f(\R^n\times \{t\})$ and $\R^n\times
%\{h(t)\}$ is at most $C$ for all $t\in \R$.  In this case, the map
%$h$ can be chosen
 % to  be   an affine
%map;
 %       to be more precise,
%there are  constants  $a\not=0$,   $b\in \R$ such that $h(t)=at+b$
%satisfies the above condition. See  Proposition 5.8 in \cite{FM}.

\b{Cor}\label{c3} {Let $A$   and $B$ be   $n\times n$ matrices whose
eigenvalues all have positive real parts.  Suppose the real part
Jordan form of $A$ is not a multiple of the identity matrix $I_n$.
 Then every quasiisometry
  $f: G_A\ra G_B$   is height-respecting.}

\end{Cor}

\b{Cor}\label{c4} {Let $A$ be a  square matrix whose eigenvalues all
have positive real parts.  If the real part Jordan form  of $A$  is
not  a multiple of the identity matrix,
  then $G_A$ is not quasiisometric to any finitely generated group.}

\end{Cor}

A group  $G$ of bijections  $g: X\ra X$ of
  a quasimetric space is a \e{unform quasim\"obius group} if there is some
  homeomorphism $\eta: [0,\i)\ra [0, \i)$ such that every element
  $g$ of $G$ is $\eta$-quasim\"obius.
%Dymarz-Peng \cite{DP}   recently proved that a   cocompact
 %  uniform
%group of  quasisimilarity  maps of $(\R^n, D_A)$ can be conjugated
%by a biLipschitz map into the
 %group of almost homotheties,  see the end of Section \ref{pro} for
  %related definitions.
   The following result follows from  Theorem
 \ref{main2}     and a  theorem  of Dymarz-Peng  \cite{DP}.
% see the end of Section \ref{pro} for
  %related definitions.

 %yields the following generalization.

\b{Cor}\label{c5} Let $A$ be a  square matrix whose eigenvalues all
have positive real parts.  Suppose that  the real part Jordan form
of $A$ is not a multiple of the identity matrix.    Let $G$ be a
uniform quasim\"obius group of $\partial G_A$ (equipped with a
visual metric).  If  the induced action of $G$ on the space of
distinct triples of $\partial G_A$  is cocompact, then $G$ can be
conjugated by a
 biLipschitz map of
   $(\R^n, D_A)$  into  the
 group of almost homotheties of  $(\R^n, D_A)$.

\end{Cor}

%\b{Cor}\label{c5} Let $A$ be a  square matrix whose eigenvalues all
%have positive real parts.  Suppose that  the real part Jordan form
%of $A$ is not a multiple of the identity matrix.    Then every
%cocompact uniform quasim\"obius group of $(\R^n, D_A)$ can be
%conjugated  by a
% biLipschitz map into  the
 %group of almost homotheties.

%\end{Cor}

When $A$ is a  Jordan block, we  describe all the quasisymmetric
maps on $(\R^n, D_A)$. Consequently, we are able to prove a
Liouville type theorem.  See Section \ref{jor} and Section
\ref{liouv}.

Theorem \ref{main2} was established in the diagonal case in
\cite{SX} and in the $2\times 2$ Jordan block case in \cite{X}. We
believe that Theorem \ref{main2}  holds true for most homogeneous
manifolds with negative curvature (HMNs), with only a few
exceptions.
 Recall  that HMNs were
 characterized by Heintze in \cite{H}:   each such manifold is
 isometric to a  solvable Lie group $S$ with a left
 invariant Riemannian metric, and the group $S$ has the form
 $S=N\rtimes \R$, where $N$ is a simply connected nilpotent Lie
 group, and the action of $\R$ on $N$ is generated by a derivation
 whose eigenvalues all have positive real parts.
  An open problem now is to establish Theorem \ref{main2} for most HMNs, and
  to construct non-biLipschitz quasisymmetric maps (of the ideal boundary)
     for the few
  exceptions.  The only  exceptions  known to the author are (those HMNs that are biLipschitz to)  the
  real and complex hyperbolic spaces:  there are
  quasisymmetric maps in the Euclidean spaces  \cite{GV}
      and the Heisenberg
  groups \cite{B}  that change Hausdorff dimensions (of certain subsets), so they can not be biLipschitz.

%There are two main questions about the rigidity of quasiisometries
%between   %finitely generated
 %solvable groups:  the quasiisometry
 %classification problem and the quasiisometric rigidity problem.
  %Farb and Mosher first solved these problems for the solvable
  %Baumslag-Solitar  groups  \cite{FM2}, \cite{FM3}
   %and more generally the so-called abelian-by-cyclic groups  \cite{FM}.
  %They found that  the rigidity  questions of (some) finitely
  %generated solvable groups is related to  those of  solvable Lie groups.
   % Then Eskin, Fisher and Whyte  \cite{EFW1}, \cite{EFW2}
    %    made remarkable progresses on the rigidity questions for
    %a large class of solvable Lie groups  and finitely generated solvable groups, including the
    %$3$-dimensional group SOL.     The solvable Lie groups
    %studied by them are not negatively curved, but they admit two transversal
    %foliations by negatively curved solvable Lie groups.
     %  Our results are the first for negatively curved solvable Lie
      % groups, which complement the above mentioned results.

\noindent {\bf{Acknowledgment}}. {I  would like to thank Bruce
Kleiner for  suggestions and stimulating  discussions. I also would
like to thank Tullia Dymarz for telling me about her joint paper
with Irine Peng \cite{DP}. Finally, I am grateful for the generous
travel support offered by the Department of Mathematical Sciences at
Georgia Southern University.}

\section{Some basic definitions}\label{prelimi}

In this section we recall some basic definitions.

    A   \e{quasimetric} $\rho$  on a set $X$
  is a function $\rho:  X\times X\ra \R$ satisfying the following
  three conditions:\newline
     (1)  $\rho(x,y)=\rho(y,x)$ for all $x, y\in
  X$;\newline
    (2)  $\rho(x, y)\ge 0$  for all $x,y\in X$,    and
    $\rho(x, y)=0$   if  and only  if  $x=y$;\newline
       (3)   there is some $M\ge 1$ such that $\rho(x, z)\le M (\rho(x,
y)+\rho(y,z))$ for all $x, y, z\in X$.\newline
   %Recall that  a quasimetric
%$\rho$  on a set $Z$ only satisfies a generalized triangle
  %inequality: there is some $M\ge 1$ such that $\rho(x, z)\le M
  %(\rho(x, y)+\rho(y,z))$ for all $x, y, z\in Z$.
       For each  $M\ge 1$,
    there is a constant $\epsilon_0>0$ such that $\rho^\epsilon$ is
      biLipschitz equivalent to a
metric for all  quasimetric $\rho$  with constant $M$  and   all
 $0<\epsilon\le \epsilon_0$,  see Proposition 14.5. in \cite{Hn}.

 For any quadruple $Q=(x,y,z,w)$ of distinct points in
    a quasimetric space $X$, the \e{cross ratio}   $\text{cr}(Q)$  of $Q$ is:
 $$\text{cr}(Q)=\frac{\rho(x, w) \rho(y,z)}{\rho(x, z) \rho(y, w)}.$$
Let $\eta: [0,\i)\ra [0,\i)$ be a homeomorphism.
 A bijection
%A homeomorphism between metric spaces
$F:X\to Y$ between two quasimetric spaces is
  \e{$\eta$-quasim\"obius} if $\text{cr}(F(Q))\le \eta(\text{cr}(Q))$ for all quadruples  $Q=(x,y,z, w)$  of
  distinct points in $X$, where $F(Q)=(F(x), F(y), F(z), F(w))$.
%Let $\eta: [0,\i)\ra [0,\i)$ be a homeomorphism.
    A bijection
%A homeomorphism between metric spaces
$F:X\to Y$ between two quasimetric spaces is
\e{$\eta$-quasisymmetric} if for all distinct triples $x,y,z\in X$,
we have
\[
   \frac{\rho(F(x), F(y))}{\rho(F(x), F(z))}\le \eta\left(\frac{\rho(x,y)}{\rho(x,z)}\right).
\]
    A map $F:  X\to Y$ is quasisymmetric if it is $\eta$-quasisymmetric
for some $\eta$.

Let $K\ge 1$ and $C>0$. A bijection $F:X\ra Y$ between two
quasimetric spaces is called a $K$-\e{quasisimilarity} (with
constant $C$) if
\[
   \frac{C}{K}\, \rho(x,y)\le \rho(F(x), F(y))\le C\,K\, \rho(x,y)
\]
for all $x,y \in X$.
   When $K=1$, we say $F$ is a \e{similarity}.
It is clear that a map is a quasisimilarity if and only if it is a
biLipschitz map. The point of using the notion of quasisimilarity is
that sometimes there is control on $K$ but not on $C$.

\section{Negatively curved $\R^n \rtimes \R$}\label{metriconb}

In this section we  first review some  basics about negatively
curved
   $\R^n\times \R$,  then
 define
the parabolic visual (quasi)metric on their  ideal boundary  and
study its  properties.  We also    recall a result of Farb-Mosher
and the main results of \cite{X} and \cite{SX}.

Let  $A$ be an  $n\times n$  matrix.  %%whose eigenvalues all
            % have positive real parts.
             Let   $\R $ act on   $\R ^n$ by
$$\R\times \R^n \rightarrow \R^n$$
             $$(t, x)\rightarrow e^{tA} x.$$
 %              \;\; (t\in \R, x\in   \R^n).$$
             We denote the corresponding semi-direct product by
             $G_A=\R^n\rtimes_A R$.  Then $G_A$ is a solvable Lie
             group.
             Recall that the group operation in $G_A$ is given by:
             $$(x_1, t_1)\cdot (x_2, t_2)=(x_1+ e^{t_1 A} x_2,  t_1
             + t_2).$$

  %Let $A$ and $G_A$ be as in the introduction.
  We will always assume that the eigenvalues of $A$ have positive real parts.
 An \e{admissible metric} on $G_A$ is a  left invariant Riemannian metric
  such that the $\R$ direction is perpendicular to the $\R^n$
  factor. The \e{standard metric} on $G_A$ is
   the left  invariant Riemannian metric
  determined  by the standard inner product on  the tangent space of
  the identity element $(0,0)\in \R^n\times \R=G_A$.
    We remark that  $G_A$ with the standard metric  does not always
       have  negative sectional curvature.
However,  Heintze's result (\cite{H}) says  that  $G_A$ has an
admissible metric with negative sectional curvature.
 Since any two  left invariant Riemannian  distances on a Lie group
 are biLipschitz equivalent,  $G_A$
with any left invariant Riemannian metric is Gromov hyperbolic.
 From now on,
   unless  indicated  otherwise, $G_A$ will always be equipped with the
 standard metric.

At a point $(x,t)\in\R^n\times\R\approx G_A$, the tangent space is
identified with $\R^n\times\R$, and the standard  metric is given by
the symmetric matrix
 \[\left(\begin{array}{cc} Q_A(t) & {\bf{0}}\\ {\bf{0}} & 1\end{array}
    \right),
  \]
where $Q_A(t)=e^{-tA^T}e^{-tA}$.  Here $T$ denotes matrix transpose.

For each $x\in \R^n$, the map $\gamma_x: \R\ra G_A$,
$\gamma_x(t)=(x,t)$  is a geodesic. We call such a geodesic a
vertical geodesic.  It can be checked that all vertical geodesics
are asymptotic as $t\ra +\infty$. Hence they define a point $\xi_0$
in the ideal boundary $\p G_A$.
    The sets $\R^n\times\{t\}$  ($t\in \R$)
are horospheres centered at $\xi_0$.
      For each  $t\in \R$, the
induced metric on  the horosphere $ \R^n\times\{t\}\subset G_A$   is
determined by the quadratic form $Q_A(t)$. This metric has distance
formula $d_{A, t}((x,t), (y,t))=| e^{-tA}(x-y)|$.  Here $|\cdot |$
denotes the Euclidean norm.
 %  The
%distance between two horospheres, corresponding to $t=t_1$ and
%$t=t_2$, is $|t_1-t_2|$. %It follows that for $(x_1,t_1),(x_2,t_2)\in
%G_A=\R^n\times\R$,
%\begin{equation}\label{eq:lowerbound}
 % d((x_1,t_1),(x_2,t_2))\ge |t_1-t_2|.
%\end{equation}

Each geodesic ray in $G_A$ is  asymptotic to either  an upward
oriented vertical geodesic or a downward oriented vertical geodesic.
The upward oriented vertical geodesics are asymptotic to $\xi_0$ and
the downward oriented vertical  geodesics are in 1-to-1
correspondence with $\R^n$. Hence $\p G_A\backslash\{\xi_0\}$ can be
naturally identified with $\R^n$.

We  next define a   parabolic visual quasimetric on
  $\p G_A\backslash\{\xi_0\}=\R^n$.
 Given $x,
y\in\R^n= \p G_A\backslash\{\xi_0\}$, the parabolic visual
quasimetric $D_A(x,y)$ is defined as follows: ${D}_A(x,y)=e^t$,
where $t$ is the smallest   real number such that at height $t$ the
two vertical geodesics $\ga_x$ and $\ga_y$ are at distance one apart
in the horosphere; that is,
  $$d_{A, t}((x,t), (y,t))=|
e^{-tA}(x-y)|=1.$$

   For each $g=(x,t)\in G_A$, the Lie group left translation $L_g$
   is an isometry of $G_A$ and fixes the point $\xi_0$.  It shifts
   all the horospheres centered at $\xi_0$  in the vertical direction by the same amount.
     It follows that the boundary map of $L_g$ is a similarity of
     $(\R^n, D_A)$.   When $g=(z,0)$,   $L_g$ leaves invariant all the horospheres centered at $\xi_0$,  and  the boundary map is
      the Euclidean
     translation by  $z$.
     %It is easy to check that
       Hence  Euclidean
translations are isometries   %and $e^{tA}$ is a similarity
    with
respect to $D_A$:
$$D_A(x+z, y+z)=D_A(x,y)\;\;\;  \text{for all}\;\;\; x,y, z\in
  \R^n.$$
    When $g=(0, t)$,  $L_g$ shifts all the horospheres centered at
    $\xi_0$ by $t$, and the boundary map is the linear
    transformation $e^{tA}$.   Hence  $e^{tA}$ is a similarity
      with similarity constant $e^t$:
     $$D_A(e^{tA}x, e^{tA}y)=e^t D_A(x,y)\;\;\; \text{for all}\;\;\; x,y\in
     \R^n  \;\;\;\text{and all} \;\;\;t\in \R.$$

 We remark that $D_A$ in general
is not a  metric, but merely  a quasimetric.
  See the  remark after the proof of Corollary \ref{hou1}.
  % If the sectional curvature of $G_A$ is $\le -1$, then $D_A$ is
  % indeed a metric.

For any integer $n\ge 2$, let
  \[J_n=\left(\begin{array}{cccccc}1 & 1 & 0 & \cdots & 0 & 0\\
  0   &  1   &  1   & \cdots   & 0   & 0\\
0   & 0  &   1    &  \cdots  &   0  &  0   \\
   &  \cdots   &   \cdots  &  \cdots &     \cdots &\\
    0 & 0 & 0 &  \cdots   &  1 & 1 \\
     0  & 0  & 0  & \cdots  & 0  & 1
     \end{array}\right)
     \]
be the $n\times n$ Jordan matrix with eigenvalue $1$.
  We write $J_n=I_n+N$. Here we omit the subscript $n$ for $N$ to simplify the notation.
   Notice that $e^{-tJ_n}=e^{-tI_n}e^{-tN}=e^{-t}e^{-tN}$.
  Hence   $D_{J_n}(x, y)=e^t$ if and only if  $t$ is the smallest
  real number satisfying
    $e^t=|e^{-tN}(y-x)|$.   For later
         use, we notice here that
 \begin{equation}\label{jormatrix}
 e^{tN}=\left(\begin{array}{cccccc}1 & t & \frac{t^2}{2!} & \cdots & \frac{t^{n-2}}{(n-2)!} & \frac{t^{n-1}}{(n-1)!}\\
  0   &  1   &  t   & \cdots   & \frac{t^{n-3}}{(n-3)!}   & \frac{t^{n-2}}{(n-2)!}\\
0   & 0  &   1    &  \cdots  &   \frac{t^{n-4}}{(n-4)!}  &  \frac{t^{n-3}}{(n-3)!}   \\
   &  \cdots   &   \cdots  &  \cdots &     \cdots &\\
    0 & 0 & 0 &  \cdots   &  1 & t \\
     0  & 0  & 0  & \cdots  & 0  & 1
     \end{array}\right).
     \end{equation}

Let $P$ be a nonsingular $n\times n$ matrix. Define a map $f: G_A
\ra G_{PAP^{-1}}$  by $f(x,t)=(Px,t)$. Then it is easy to check that
$f$ is a  Lie group isomorphism.  Hence $f$ is an isometry if
$G_{PAP^{-1}}$ is equipped with the standard metric and $G_A$ has
the admissible metric in which $P^{-1}e_1, \cdots, P^{-1}e_n,\,
e_{n+1}$ is orthonormal
 at the identity element of $G_A$.   Here $e_1, \cdots, e_n$ denote
 the standard basis of $\R^n$, and  $e_{n+1}$ is the standard basis for $\R$.
Hence, $G_A$ with any admissible metric is isometric to
  $G_{PAP^{-1}}$  with the standard metric  for some nonsingular matrix $P$.  By Heintze's result \cite{H}, there is
  a nonsingular matrix $P$ such that
$G_{PAP^{-1}}$  with the standard metric  has negative sectional
curvature.

 Now we suppose both $G_A$ and  $G_{PAP^{-1}}$  are equipped with the standard
 metric.  Then it is easy to check that
 for each $t\in \R$, the restricted map
   $$f|_{\R^n\times \{t\}}:
 (\R^n\times \{t\},\,   d_{A, t})  \ra (\R^n\times \{t\},\,  d_{PAP^{-1}, t})$$
  is   $K$-biLipschitz,
     %where the
   %domain $\R^n\times \{t\}$ is equipped with the Riemannian metric
  %induced from $G_A$ and  the target $\R^n\times \{t\}$   is equipped with the Riemannian metric
  %induced from   $G_{PAP^{-1}}$.
   %The biLipschitz constant can be taken to be
     where $K:=\max\{||P||,
   ||P^{-1}||\}$.    Here    $||M||=\sup\{|Mx|:   x\in \R^n, |x|=1\}$
       denotes the operator norm  of
    an $n\times n$ matrix $M$.
    %$||M||=\sup\{|Mx|:   x\in \R^n, |x|=1\}$.
We next recall a more general result by  Farb-Mosher \cite{FM}.

  %Given a $n\times n$ matrix  $A$, the real part Jordan form    of
  %$A$  is obtained from the Jordan form of $A$ by replacing each
  %diagonal entry with  its real part  and reordering  to make it
  %canonical.

\b{Prop}\label{farb-m} \e{(Proposition 4.1 in \cite{FM})}
   Let $A$ and $B$ be two $n\times n$ matrices. Suppose there are  constants
   $r, s>0$ such that $rA$ and $sB$ have the same real part Jordan
   form.  Then there is a height-respecting quasiisometry
    $f: G_A\ra G_B$.  To be more precise,  there exist an $n\times n$ matrix
    $M$ and $K\ge 1$ such that for each $t\in \R$, the map
     $v\ra Mv$  is a $K$-biLipschitz  homeomorphism from
      $(\R^n, d_{A, t})$ to $(\R^n, d_{B, \frac{s}{r}t})$;  it
      follows that the map  $f: G_A\ra G_B$  given by
      $$(x,t)\longmapsto \left (Mx, \frac{s}{r}\cdot t\right)$$
      is biLipschitz with biLipschitz constant
        $\sup\{K, \frac{s}{r}, \frac{r}{s}\}$.

    \end{Prop}

\b{Cor}\label{hou1} Suppose we are in the setting of  Proposition
\ref{farb-m}. Assume further that $r=1$ and $G_A$ has negative
sectional curvature. Then:\newline
     (1) the boundary map $\p f: (\R^n, D^s_A)\ra  (\R^n, D_B)$ is
     biLipschitz;\newline
(2)  $f$ is an almost similarity.

\end{Cor}

\b{proof} (1) We observe that the boundary map is given by
  $\p f(x)=M x$.   Let
$x,y\in\R^n$   and   assume $D^s_A(x,y)=e^t$. Then
$D_A(x,y)=e^{t/s}$. By the definition of $D_A$, we have $d_{A,
t/s}((x,t/s), (y,t/s))=1$. Since $G_A$ has pinched negative
sectional curvature, there is a constant $a$ depending only on the
curvature bounds of $G_A$, such that
 $d_{A, t'}((x,t'), (y,t'))<1/K$  for $t'>t/s+a$  and
$d_{A, t'}((x,t'), (y,t'))>K$  for $t'<t/s-a$. It  now follows from
Proposition \ref{farb-m}
  that $d_{B, t''}((Mx, t''), (My, t''))<1$
    for $t''>t+sa$ and $d_{B, t''}((Mx, t''), (My, t''))>1$
  for  $t''< t-sa$.  By  the definition of $D_B$ we have
    $e^{-sa}e^t\le D_B(Mx, My)\le e^{sa}e^t$.
      Hence
$\p f: (\R^n, D^s_A)\ra  (\R^n, D_B)$ is
     biLipschitz  with biLipschitz constant $e^{sa}$.

  (2) Let $p=(x_1, t_1)$, $q=(x_2, t_2)\in G_A$ be arbitrary.
   We may assume $t_1\le t_2$.  If $x_1=x_2$, then it is clear from
   the definition of $f$ that $d(f(p), f(q))=s\cdot d(p,q)$. So we assume
    $x_1\not=x_2$ and that  %$D_A(x_1, x_2)=e^{t_0}$ for some $t_0$.
       $d_{A, t_0}((x_1, t_0), (x_2, t_0))=1$ for some $t_0$.
First assume $t_0\le t_2$.  Then $d((x_1, t_2), q)< d_{A, t_2}((x_1,
t_2),  q)\le 1$ as $G_A$ has negative sectional curvature.
 By the triangle inequality, we have $|d(p,q)-(t_2-t_1)|\le 1$.
  By Proposition \ref{farb-m}, $d((M x_1, st_2), f(q))\le  d_{B, st_2}((Mx_1,
  st_2), f(q))\le K$.
 By the triangle inequality again we have
   % $|d((Mx_1, st_1), (Mx_2, st_2))-(st_2-st_1)|\le K$,
    %that is,
       $|d(f(p), f(q))-(st_2-st_1)|\le K$. Hence
    $|d(f(p), f(q))-s\cdot  d(p, q)|\le s+K$.

       Next we assume $t_0>t_2$.  By Lemma 6.3 (1) of \cite{SX} we
       have  $|d(p, q)-(t_0-t_1)-(t_0-t_2)|\le C_1$ for some constant $C_1$
       depending only on the curvature bounds of $G_A$.
  By Lemma 6.2 of \cite{SX}, the point $(x_1, t_0)$ is a
  $C_2$-quasicenter of $x_1, x_2, \xi_0\in \p G_A$ for some constant
  $C_2$ depending only on the curvature bounds of $G_A$.
 Since $f$ is a quasiisometry, $f(x_1, t_0)=(M x_1, st_0)$ is a $C_3$-quasicenter
 of $M x_1, M x_2, \eta_0\in \p G_B$ (here $\eta_0$ denotes the point in $\p G_B$ corresponding to
    upward oriented vertical geodesics), where $C_3$ depends
 only on $C_2$, the quasiisometry constants of $f$ and the Gromov
 hyperbolicity constant of $G_B$.
  Similarly, the point $(M x_2, st_0)$ is also a
  $C_3$-quasicenter of
 $Mx_1, Mx_2, \eta_0\in \p G_B$. Now consider the geodesic triangle
 consisting of
 $\{M x_1\}\times \R$, $\{M x_2\}\times \R$ and  a geodesic joining $M x_1$, $M
 x_2$.
 Notice that $f(p)\in \{M x_1\}\times \R$ lies between $Mx_1$ and $(Mx_1, st_0)$
 and $f(q)\in \{M x_2\}\times \R$ lies  between $Mx_2$ and  $(Mx_2, st_0)$.  It
 follows that
 \b{align*}
  \;  & |d(f(p), f(q))-(st_0-st_1)-(st_0-st_2)|\\
  &=|d(f(p), f(q))-d(f(p), (Mx_1, st_0))-d(f(q), (M x_2, st_0))| \le C_4
 \end{align*}
  for some constant $C_4$ depending only on $C_3$
 and the Gromov hyperbolicity constant of $G_B$.
This combined with
   $|d(p, q)-(t_0-t_1)-(t_0-t_2)|\le C_1$
 implies $|d(f(p), f(q))-s\cdot  d(p,q)|\le C_4+s
C_1$.

\end{proof}

We notice that Corollary  \ref{hou1}  (1) implies that  $D_A$ is
  indeed a
quasimetric:  by Heintze's result, there is some nonsingular $P$
such that $G_{PAP^{-1}}$ has pinched negative sectional curvature
and hence $D_{PAP^{-1}}$ is a quasimetric (this can be proved by the
arguments in \cite[p.124]{CDP} or by using the relation between
parabolic visual quasimetric and visual quasimetric \cite[section
5]{SX}); since
 $(\R^n,  D_A)$ and $(\R^n,  D_{PAP^{-1}})$ are biLipschitz, $D_A$ is also a
 quasimetric.

%We notice that a quasiisometry $f: G_A\ra G_B$ is height-respecting
%if and only if there is some $\epsilon>0$ such that  the boundary
%map $\p f: (\R^n, D_A)\ra (\R^n, D_B^\epsilon)$
 %is biLipschitz.   The proof
  % of this statement is similar to that  of Lemma 6.4  in \cite{SX} (the notion of height-respecting is more restrictive there).

   %  Because of Proposition \ref{farb-m}, we will assume from now on
    % that   $A$ is already in its real part  Jordan form  for the matrix $A$ in $G_A$.

Let $A$ be an  $n\times n$ matrix in real part Jordan form with
positive eigenvalues
 $$\la_1<\cdots<\la_{k_A}.$$
  Let $V_i\subset\R^n$ be the generalized eigenspace of $\la_i$, and
let $d_i=\dim V_i$.

   If $k:=k_A\ge 2$, we write $A$ in the block diagonal form
    $A=[A_1, \cdots, A_{k}]$,  where $A_i$ is the block
    corresponding to the eigenvalue  $\la_i$;
 we also denote $A'=[A_1, \cdots, A_{k-1}]$.  % and $B'=[B_1, \cdots, B_{k_B-1}]$.
  Correspondingly,  $\R^n$ admits the decomposition  $\R^n=V_1\times \cdots \times V_k$. Hence  each point $x\in \R^n$ can be written
$x=(x_1, \cdots, x_k)$, where $x_i\in V_i$.
   Observe that, for each $x_k\in V_k$, if we identify $V_1\times
   \cdots \times V_{k-1}\times \{x_k\}$  with  $V_1\times
   \cdots \times V_{k-1}$, then the restriction of $D_A$ to
 $V_1\times
   \cdots \times V_{k-1}\times \{x_k\}$  agrees with $D_{A'}$.
  It is not hard to   check  that for all $x_k, y_k\in V_k$, the following holds for  the Hausdorff
  distance   with respect  to
the quasimetric $ D_A$:
\begin{equation}\label{eq:1}
  HD(V_1\times \cdots \times V_{k-1}\times \{x_k\}, V_1\times \cdots \times V_{k-1}\times \{y_k\})=D_{A_k}(x_k,  y_k).
\end{equation} %=(d_Y(y_1, y_2))^{\frac{\alpha_2}{\alpha_1}}.$$
Also, for any $x=(x_1, \cdots,  x_k)\in \R^{n}$ and any $y_k\in
V_k$,
\begin{equation}\label{eq:2}
    D_A(x, V_1\times \cdots \times V_{k-1}\times \{y_k\})=D_{A_k}(x_k,y_k).
\end{equation}

  When $k=1$, that is,  when  $A$ has only one eigenvalue $\la:=\la_1>0$, the
  matrix $A$ also has  a block diagonal  form
    $A=[\la I_{n_0}, \la I_{n_1}+N,\cdots, \la I_{n_r}+N]$, where $n_0\ge 0$  and $\la I_{n_i}+N$ is a Jordan
 block.  We allow the case $A=\la I_n$.
  We write  a point $p\in \R^n$ as $p=(z,(x_1, y_1), \cdots, (x_r,
  y_ r))^T$, where $T$ denotes  matrix transpose, $z\in \R^{n_0}$
  corresponds to $\la I_{n_0}$   and  $(x_i, y_i)^T\in \R^{n_i}$
  ($x_i\in \R^{{n_i}-1}$, $y_i\in \R$)
  corresponds to $\la I_{n_i}+N$.      %and    $\{(x_i, y_i): x_i=0\}$ is
     %the eigenspace of $\alpha I_{n_i}+N$.
    Set
  $$\R^{n_0 +r}=\{p=(z,(x_1, y_1), \cdots, (x_r,
  y_ r))^T\in \R^n:  x_1=\cdots=x_r={\bf{0}}\},$$
    and let $\pi_A:  \R^n\ra \R^{n_0+r}$ be the  projection  given
    by:
    $$\pi_A(p)=(z, ({\bf{0}}, y_1), \cdots, ({\bf{0}},  y_r))^T   \;\, \text{for} \;\,p=(z,(x_1, y_1), \cdots, (x_r,
  y_ r))^T\in \R^n.$$
      %Notice that  $\R^{n_0 +r}$ is the eigenspace of $A$.
         Set
           $$A(1)=[\la I_{n_1-1}+N,\cdots, \la I_{n_r-1}+N],$$
            where $\la I_1+N$ is understood to be $\la I_1$.

\b{Le}\label{l11} {
  The restriction of $D_A$ to the fibers of $\pi_A$ agrees with
    $D_{A(1)}$. To be more precise,
   for all $p=(z,(x_1, y_1), \cdots, (x_r,
  y_ r))^T$, $p'=(z,(x'_1, y_1), \cdots, (x'_r,
  y_ r))^T$  we have
    $$D_A(p, p')=D_{A(1)}(x, x'),$$
     where $x=(x_1, \cdots, x_r)^T$ and $x'=(x'_1, \cdots, x'_r)^T$.

}

\end{Le}

\b{proof} Assume $D_A(p, p')=e^t$ and $D_{A(1)}(x, x')=e^s$.
  By the definition,
  $s$ is the smallest real number such that $|e^{-sA(1)}(x'-x)|=1$.
 We calculate
$$e^{-sA(1)}(x'-x)=e^{-\la s} (e^{-sN_{n_1-1}}(x'_1-x_1), \cdots,e^{-sN_{n_r-1}}(x'_r-x_r))^T.$$
  Similarly,
  $t$ is the smallest real  number such
  that $|e^{-tA}(p'-p)|=1$.  We calculate
  $$e^{-tA}(p'-p)=e^{-\la t}({\bf{0}}, (e^{-tN_{n_1-1}}(x'_1-x_1),
  0), \cdots, (e^{-tN_{n_r-1}}(x'_r-x_r),
  0))^T.$$
    It follows that
the two equations  $|e^{-sA(1)}(x'-x)|=1$  and
  $|e^{-tA}(p'-p)|=1$
 are the same. Hence
$s=t$.

\end{proof}

\b{Le}\label{l12} {The following hold for all $y, y'\in
\R^{n_0+r}$:\newline (1)  the Hausdorff distance
  $HD_{D_A}(\pi_A^{-1}(y), \pi_A^{-1}(y'))=|y-y'|^{\frac{1}{\la}}$; \newline
     (2) for any $p\in \pi_A^{-1}(y)$,
     we have
      $D_A(p,\pi_A^{-1}(y'))=|y-y'|^{\frac{1}{\la}}$.

}

\end{Le}

\b{proof} Let $p=(z,(x_1, y_1), \cdots, (x_r,
  y_ r))^T\in \pi_A^{-1}(y)$, $p'=(z',(x'_1, y'_1), \cdots, (x'_r,
  y'_ r))^T\in  \pi_A^{-1}(y')$, where $y$ and $y'$ are written
  $y=(z, y_1, \cdots, y_r)$, $y'=(z', y'_1, \cdots, y'_r)$.
    Assume  $D_A(p, p')=e^t$. Then $t$ is the smallest  real number
    such that
$$\big\vert\left(z'-z, e^{-tN_{n_1}}(x'_1-x_1, y'_1-y_1)^T, \cdots,
e^{-tN_{n_r}}(x'_r-x_r, y'_r-y_r)^T\right)\big \vert=e^{\la t}.$$
      Notice  that  the
last entry of   $e^{-tN_{n_i}}(x'_i-x_i, y'_i-y_i)^T$
   is  $y'_i-y_i$, which is independent of $t$. It follows that
    $e^{\la t}\ge |(z'-z, y'_1-y_1, \cdots, y'_r-y_r)|=|y'-y|$,
    and hence $D_A(p, p')=e^t\ge |y'-y|^{\frac{1}{\la}}$.

Set $t_0=\ln |y'-y|/\la$. Then $e^{\la t_0}=|y'-y|$.
  Now let $p=(z,(x_1, y_1), \cdots, (x_r,
  y_ r))^T\in \pi_A^{-1}(y)$  be arbitrary.
    Since the matrix  $e^{-t_0N_{n_i}}$ is nonsingular, the equation
 $$e^{-t_0N_{n_i}}(u_i, v_i)^T=(0, \cdots, 0, y'_i-y_i)^T$$
 has a unique
 solution  $(u_i, v_i)^T$, where $u_i\in \R^{n_i-1}$ and $v_i\in \R$. Notice that
  $v_i=y'_i-y_i$.    Set $x'_i=u_i+x_i$ and
   $p'=(z',(x'_1, y'_1), \cdots, (x'_r,
  y'_ r))^T$. Then $p'\in\pi_A^{-1}(y')$  and
  $$e^{-t_0 A}(p'-p)=e^{-t_0 \la}\left(z'-z, ({\bf{0}}, y'_1-y_1), \cdots, ({\bf{0}}, y'_r-y_r)\right)^T.$$
  It follows that $t_0$ is a  solution of
  $|e^{-tA}(p'-p)|=1$ and so
   $D_A(p, p')\le e^{t_0}=|y-y'|^{\frac{1}{\la}}$. This together with the first paragraph implies
       $D_A(p, p')=|y-y'|^{\frac{1}{\la}}$.
  So each point $p\in \pi_A^{-1}(y)$ is within
  $|y-y'|^{\frac{1}{\la}}$  of
  $\pi_A^{-1}(y')$.  Similarly, every point $p'\in \pi_A^{-1}(y')$ is
  also within   $|y-y'|^{\frac{1}{\la}}$  of
  $\pi_A^{-1}(y)$.  Therefore, (1) holds.

   (2) also follows from the above two paragraphs.

\end{proof}

The following two results will be used in the proof of Theorems
\ref{main1} and \ref{main2}.   They are the basic steps in the
induction.

\begin{theorem} \e{(Theorem 4.1 in  \cite{SX})} \label{sx}
 Suppose   $A$   is  diagonal
 with positive eigenvalues  $\alpha_1<\alpha_2<\cdots <\alpha_r$  \e{($r\ge 2$)}.  % Suppose $A$ is not
  % a multiple of the identity matrix.
   Then
   every $\eta$-quasisymmetry  $F:  (\R^n, D_A)\ra (\R^n, D_A)$
        is  a  $K$-quasisimilarity, where $K$ depends only on $\eta$
        and $r$.
  % $\partial G_A\setminus\{\xi_0\}$
 %       with respect to the parabolic visual metric
%$D_e$.  % and preserves horizontal foliations.
\end{theorem}

\b{theorem} \e{(Theorems 4.1 and 5.1 in    \cite{X})} \label{xthm}
 {%Let
  %\[A=\left(\begin{array}{cc} 1 & 1\\ 0 & 1\end{array}\right).
  %\]
   Every  $\eta$-quasisymmetric map $F: (\R^2, D_{J_2})\ra (\R^2, D_{J_2})$  is a
 $K$-quasisimilarity, where $K$ depends only on $\eta$.
   Furthermore, a bijection $F:(\R^2, D_{J_2})\ra (\R^2, D_{J_2})$
   is a quasisymmetric map   if  and only of it
has the following form: $F(x,y)=(ax+c(y), ay+b)$
  for all $(x,y)\in \R^2$, where $a\not=0$,  $b$ are constants and
  $c: \R\ra \R$ is a Lipschitz map.}
\end{theorem}

\section{$Q$-variation on the ideal boundary}\label{qvar}

In this section we introduce the main tool in the proof of the main
results: $Q$-variation   for maps between quasimetric spaces.
 It is a discrete version of the notion of capacity.  The advantage
 of this notion is that it makes sense for quasimetric spaces and
  does not require the existence of rectifiable curves. We  remark
  that, while dealing with ideal boundary of negatively curved
  spaces,  very often    either one has to work with quasimetric spaces in which the
  triangle inequality is not available, or   one needs to work with
  metric spaces that contain no rectifiable curves.   Both
  scenarios  are unpleasant from the point of view of classical
  quasiconformal analysis.
 % As a tool, $Q$-variation could be compared with Pansu's  modulus  \cite{P}, but
%seems slightly easier to work with in our  context.

The notion of  $Q$-variation  is due to Bruce Kleiner \cite{K}.

Let $(X, \rho)$ be a quasimetric space and $L\ge 1$.
 A subset $A\subset X$ is  called an $L$-\e{quasi-ball} if there is
 some $x\in X$ and some $r>0$ such that $B(x,r)\subset A\subset B(x,
 Lr)$.    Here $B(x,r):=\{y\in X:   \rho(y,x)<r\}$.  % denotes the \lq\lq open ball"  with center $x$ and radius $r$.}

 For  any ball $B:=B(x,r)$ and any $\kappa>0$,  we  sometimes  denote $B(x,
 \kappa r)$ by  $\kappa B$.

 For a subset $E$ of a quasimetric space $(Y, \rho)$,
   the $\rho$-diameter of $E$ is
   $$\text{diam}_\rho(E):=\sup \{\rho(e_1, e_2): e_1, e_2\in E \}.$$
 Let $u:  (X, \rho_1)\ra (Y, \rho_2)$ be a map between two quasimetric spaces.
   For any subset $A\subset X$, the oscillation of $u$ over $A$ is
     $$\text{osc}(u|_A)=\text{diam}_{\rho_2}(u(A)).$$
   Let $Q\ge 1$.
For
     a     collection of  disjoint subsets   $\mathcal{A}=\{A_i\}$   of
     $X$, the  \e{$Q$-variation of
   $u$  over ${\mathcal{A}}$},   denoted  by  $V_Q(u,{\mathcal{A}})$,   is the quantity
$$\sum_i[osc(u|_{A_i})]^Q.$$
 For $\delta>0$ and $K\ge 1$,  set
  %the \e{$(Q,K)$-variation $V^\delta_{Q,K}(u)$   of $u$}  is
     $$V^\delta_{Q,K}(u)=\sup\{V_Q(u, {\mathcal{A}})\}, $$
      where ${\mathcal{A}}$
       ranges over  all disjoint  collections  of $K$-quasi-balls in
       $(X,\rho_1)$
         with $\rho_1$-diameter at most $\delta$.
           Finally,  the \e{$(Q,K)$-variation $V_{Q,K}(u)$   of $u$}  is
             $$V_{Q,K}(u)=\lim_{\delta\ra 0}V^\delta_{Q,K}(u).$$
      %Here the $D$-diameter of a subset $A\subset X$ is
      %$diam_D(A):=\sup \{D(a_1, a_2): a_1, a_2\in A \}$.

We notice that $V_{Q,K}(u|_{E_1})\le V_{Q,K}(u|_{E_2})$  whenever
 $E_1\subset E_2\subset X$.

There are useful variants of this definition, for instance one can
look at the infimum over all coverings. Or one can take the infimum
over all coverings followed by the sup as the mesh size tends to
zero. As a tool, $Q$-variation could be compared with Pansu's
modulus  \cite{P}, but seems slightly easier to work with in our
context.
% The definition preforms the same function as Pansu's modulus
 %   \cite{P}, but it seems  easier to digest.

Since quasisymmetric maps send quasi-balls to quasi-balls
quantitatively,  it is easy to see that $Q$-variation
  is a  quasisymmetric   invariant.   To be more precise, we
  recall

\b{Le}\label{l3.0}  \e{(Lemma  3.1  in \cite{X})}
   {Let  $X$ be a bounded quasimetric space  and  $F:   X\ra Z$  an $\eta$-quasisymmetric map.
          Then  for every map  $u:  X\ra Y$ we
   have
    $V_{Q, K}(u)\le V_{Q, \eta(K)}(u\circ F^{-1})$.

   }

   \end{Le}

  We next calculate the $Q$-variation of certain functions defined on the
  ideal boundary of negatively curved $\R^n\rtimes \R$. These
  calculations will be used in the next section to show that certain
  foliations on the ideal boundary are preserved by quasisymmetric
  maps.

  For later use   we recall that,  for any
  $Q>1$, any
  integer $k\ge 1$ and any nonnegative numbers $a_1, \cdots, a_k$,
    Jensen's  inequality states
   $$\frac{\sum_{i=1}^{k} a_i^Q}{k}\ge
   \left(\frac{\sum_{i=1}^{k}a_i}{k}\right)^Q,$$
    and equality holds if and only if  all the
    $a_i$'s   are
  equal.
     In our applications, the $a_i's$ will   be the oscillations of a function $u$
  along a ``stack" of quasi-balls.

Let $A$ be an $n\times n$ matrix in real part Jordan form with
positive eigenvalues
$$ \la_1<\la_2<\ldots <\la_k\,,
$$
 let $V_i\subset\R^n$ be the generalized eigenspace of $\la_i$,
and let $d_i=\dim V_i$.  Then $\R^n$ admits the decomposition:
   $\R^n=V_1\times \cdots \times V_k$.
     Since
$e^{tA}$ is a linear transformation with
$\text{det}(e^{tA})=e^{t(\sum_i\,d_i\la_i)}$,
 for any subset $U\subset \R^n$, we have
 $\vol(e^{tA}(U))=e^{t(\sum_i\,d_i\la_i)}\vol(U)$.

There are constants   $C_1, C_2, C_3$
  depending only on the dimension $n$ with the following properties.
    If $B:=B(o,1)\subset\R^n$ is the unit ball (in the Euclidean metric), and $t\leq
-1$, then
\begin{equation} \label{eqn-inoutradius}
B(o,   C_1e^{t\la_k}|t|^{-n+1})\subset e^{tA}\,B\subset B(o,
C_2e^{t\la_1}|t|^{n-1})\,,
\end{equation}
while
    \b{equation} \label{eqn-volume}
\vol(e^{tA}B)=C_3\,e^{t(\sum_i\,d_i\la_i)}\,.
\end{equation}

Let  $S=\prod_{i=1}^n [0,1]\subset \R^n$  be the unit cube.
 We notice that  both $S$  and  $B$   are  $K_0$-quasi-balls with respect to $D_A$  for
 some $K_0$ depending only on  $A$. Hence there is some $r>0$ such
 that
  $B_A(o, r)\subset B\subset B_A(o, K_0r)$. Here  the subscript $A$
  refers to $D_A$.  Also recall that $D_A$ is a quasimetric:  there
  is a constant $M\ge 1$ such that $D_A(x,z)\le
  M(D_A(x,y)+D_A(y,z))$ for all $x,y, z\in \R^n$.

In the following, when we say a subset $E\subset \R^n$ is convex, we
 mean it is convex with respect to the Euclidean metric.
  The  continuity of a function $u: E\ra \R$ is with respect to the
  topology induced from the usual topology on $\R^n$.

\b{Le}\label{l3.1} %%There is a constant
  %$K_0\ge 1$ depending only on the matrix $A$ with the following
    %property.
       Let  $E\subset \R^n$ be a convex open subset.
  If  $u:  (E,  D_A)  \ra\R$ is a nonconstant continuous function,    then $V_{Q,
K}(u)=\i$ for all $Q<\frac{\sum_i\,d_i\la_i}{\la_k}$  and   all
$K\ge K_0$.
   %  and all subsets $E\subset \R^n$  containing a sufficiently
  %large Euclidean ball about the origin.

\end{Le}

\b{proof} %Suppose $u:\R^n\ra\R$ is a nonconstant continuous
   %function,
    %Let  $E_0\subset \R^n$ be a large Euclidean ball containing  two points
   Let  $p, q\in E$   with
      $u(p)\neq u(q)$. Let $C\subset E$ be a fixed
  cylinder   with axis
    $\ol{pq}$, such that the minimum of $u$ on one cap of $C$ is
strictly greater than its maximum on the other cap. We pack $C$ with
translates of $e^{tA}B$, for $t\ll 0$, as follows.  First pick
  a maximal set of  lines
   $\mathcal{L}=\{L_j\}$  in $\R^n$
  satisfying the following  conditions: \newline
   (1) each line is   parallel to
    $\ol{pq}$;\newline
(2) each line  intersects $C$;\newline
      (3)  the Hausdorff distance (with respect to
    $D_A$) between any two of the lines  is  at least $2MK_0re^t$.\newline
     The maximality implies that for each $x\in C$,  we have $D_A(x, L_j)\le
  2MK_0re^t$ for some $j$.
   For each $j$,  consider a translate  $B_j$ of $e^{tA}B$ centered at
   some point on $L_j$.  Then we move $B_j$ along $L_j$  (in both directions)  by translations
      until the translates just touch $B_j$. Repeat this and we
      obtain a ``stack" of  $K_0$-quasi-balls  centered on $L_j$.
 Do this for each $j$ and we obtain  a packing   $\mathcal{P}=\{P\}$   of $C$ by translates of
 $e^{tA}B$, after removing those that are disjoint from $C$.

We claim that
 the collection $\mathcal{P}$    covers
 a fixed fraction of the volume
of $C$.  To see this, first notice that the $D_A$-distance between
the centers  $x_1$, $x_2$  of two consecutive $K_0$-quasi-balls
along $L_j$ is at most $M(K_0re^t+K_0re^t)=2MK_0re^t$, due to the
generalized triangle inequality for $D_A$.   Assume  $D_A(x_1,
x_2)=e^s$.
   Then
     $$e^{(\ln r-s)A}(x_2-x_1)\in e^{(\ln r-s)A}\ol{B}_A(o, e^s)=\ol{B}_A(o,r)\subset \ol{B}\subset
 \ol{B}_A(o, K_0 r).$$
     Since $\ol{B}$ is convex,  the line segment  joining $o$ and
$e^{(\ln
 r-s)A}(x_2-x_1)$
 is contained in  $\ol{B}\subset
 \ol{B}_A(o, K_0 r)$.
 It follows that   the segment joining $o$ and $x_2-x_1$ lies in  $ e^{(s-\ln r)A}
 \ol{B}_A(o, K_0 r)=\ol{B}_A(o,
 K_0 e^s)$.  Hence $\ol{x_1x_2}\subset \ol{B}_A(x_1, K_0 e^s)$.  This shows that every
 point on $L_j\cap C$ is within $K_0e^s\le K_1:=2MK^2_0re^t$ of  the center of
 some $P\in \mathcal{P}$.   Now the choice of the lines $\{L_j\}$
 and the generalized triangle inequality for $D_A$ imply that
  $C$ is covered by $D_A$-balls with radius $K_2:=M(2MK_0re^t+K_1)$
  and centers the centers of $\{P\}$.    Since the volumes of
  $e^{tA}B$ and  $B_A(o,  K_2)$ are comparable, the claim follows.

   The number of
$K_0$-quasi-balls in $\mathcal{P}$ along each line $L_j$   is
$\lesssim e^{-t\la_k}|t|^{n-1}$ in view of the estimate
(\ref{eqn-inoutradius}).  By Jensen's inequality, the
  $Q$-variation of  $u$ for   the $K_0$-quasi-balls along
  $L_j$ is at least as large as the $Q$-variation   when  the oscillations of $u$ on these
  quasi-balls are equal.  This common oscillation is $\gtrsim
  e^{t\la_k}|t|^{-n+1}$.
   Since $\mathcal{P}$ covers a fixed fraction of $C$, the
    cardinality of  $\mathcal{P}$ is  $\gtrsim e^{-t(\sum_i\,d_i\la_i)}$.
       Hence the
$Q$-variation of $u$ on $\mathcal{P}$   is
$$
\gtrsim e^{-t(\sum_i\,d_i\la_i)}\left(e^{t\la_k}|t|^{-n+1} \right)^Q
$$
$$
=e^{t(Q\la_k-\sum_i\,d_i\la_i)}|t|^{(-n+1)Q}\,,
$$
which $\ra\infty$ as $t\ra -\infty$ for  $Q<
\frac{\sum_i\,d_i\la_i}{\la_k}$.    Hence $V_{Q, K}(u)=\i$.
  %  for any subset $E\subset \R^n$ containing the cylinder $C$.

\end{proof}

 Notice that $\frac{\sum_i\,d_i\la_i}{\la_k}<n$ if $k\ge 2$  and
 $\frac{\sum_i\,d_i\la_i}{\la_k}=n$  if $k=1$. Hence we have the following:

\b{Cor}\label{c3.1} {Suppose $k=1$.
  Let $E\subset \R^n$ be a convex open subset.
 If  $u:(E, D_A) \ra\R$ is a
nonconstant continuous function,    then $V_{Q, K}(u)=\i$ for all
$Q<n$  and   all $K\ge K_0$.
%  and all subsets $E\subset \R^n$ containing
  %a sufficiently large Euclidean ball about the origin.

}

\end{Cor}

\b{Le}\label{l3.2} Let $E\subset \R^n$ be a convex open subset.
    Let
$u:   (E,  D_A) \ra\R$ be a continuous function.
 Suppose  there is  an affine subspace  $W$
 parallel to the subspace $\prod_{i\le l}V_i$ such that $u\vert_{W\cap E}$ is
 not constant.
Then $V_{Q, K}(u)=\i$ for all $Q<\frac{\sum_i\,d_i\la_i}{\la_l}$ and
all $K\ge K_0$.
    %  and all subsets $E\subset \R^n$  containing a
  %sufficiently large Euclidean ball about the origin.

\end{Le}

\b{proof}

Note that in the proof of Lemma \ref{l3.1},
  if $\ol{pq}$ is parallel
to the subspace $\prod_{i\le
 l}V_i$, then the number of quasi-balls in
 ${\mathcal{P}}$    along a line $L_j$  is $\lesssim
e^{-t\la_l}|t|^{n-1}$, so the lower bound on $Q$-variation becomes
%\begin{equation}
%\label{eqn-jlowerbound}
\[
C\,e^{t(Q\la_l-\sum_i\,d_i\la_i)}|t|^{(-n+1)Q}\,,
\]
%\end{equation}
which tends to $\infty$ as $t\ra-\infty$  if   $Q<
\frac{\sum_i\,d_i\la_i}{\la_l}$.

\end{proof}

  Let $\pi:  \R^n=V_1\times \cdots \times V_k\ra V_k$ be the natural
  projection.

\b{Le}\label{l3.3}
  %Let $\pi: (\R^n, D_A)\ra V_k$ be the projection onto $V_k$, and
    Let    $\pi': V_k\ra \R$ be  a coordinate function on $V_k$,  and  set
  $u=\pi'\circ \pi$.
  %Suppose $u:  (\R^n,  D_A)  \ra\R$ is a linear function factoring
%through the projection to $V_k$.
   Then $V_{Q, K}(u|_E)=0$ for all
$Q>\frac{\sum_i\,d_i\la_i}{\la_k}$,   all $K\ge
  K_0$  and all bounded subsets    $E\subset \R^n$.

\end{Le}

\b{proof} Let $E$ be a  bounded open subset.
 Let
$0<\delta<<1$ and
  %Since  $u:\R^n\ra\R$ is a linear function factoring through the
%projection to $V_k$,  for any $t<-1$, if
   $\{B_j\}_{j\in I}$  be  a packing of $E$ by    $K$-quasi-balls  with size $<\delta$.
  Then for each $j$ there is some  $x_j\in \R^n$ and some $t_j$ such that
  $$B_A(x_j, e^{t_j})\subset  B_j\subset B_A(x_j, Ke^{t_j}).$$
 Since $B_A(o, r)\subset B\subset B_A(o,K_0 r)$, we   have
 $e^{t'_j A}B\subset B_A(o, e^{t_j})$ and $B_A(o, Ke^{t_j})\subset
 e^{t''_j A}B$, where $t'_j=t_j-\ln r-\ln K_0$ and $t''_j=t_j-\ln
 r+\ln K$.  Set  $B'_j=x_j+e^{t'_j A}B$  and
   $B''_j=x_j+e^{t''_j A}B$.
   Then
   $B'_j\subset B_j\subset B''_j$.   It follows that
    %the oscillation of $u$ on $B_j$ is
$${\text{osc}}(u|_{B_j})\le \text{osc}(u|_{B''_j})
\lesssim e^{t''_j\la_k}|t''_j|^{d_k-1},
$$
  and    %the contribution of $u$  on $B_j$ to the $Q$-variation is
$$({\text{osc}}(u|_{B_j}))^Q
\lesssim e^{t''_j(Q\la_k)}|t''_j|^{Q(d_k-1)}\,\lesssim
e^{t'_j(Q\la_k)}|t'_j|^{Q(d_k-1)}\,.
$$
If $Q>\frac{\sum_i d_i\la_i}{\la_k}$,
    then this will be $\lesssim
(\vol(B'_j))^s\le   (\vol(B_j))^s$ for  $s=\frac{Q\la_k+\sum_i
d_i\la_i}{2\sum_i d_i\la_i}>1$, which implies that the $Q$-variation
  is  zero.

\end{proof}

%Putting all this together, one gets that the collection of
%continuous functions with locally finite $p$-variation for all
%$p>\frac{\sum_id_i\la}{\la_k}$ detects the affine subspaces parallel
%to  $\sum_{i<k}\,V_i$.

For the rest of this section, we will assume $k=1$  and use the
notation introduced before Lemma \ref{l11}.

\b{Le}\label{l3.4}
  %There is a constant $K_0\ge 1$  with the following property.
  Suppose $A$ has only one eigenvalue $\la>0$.
  % Suppose $r=1$.
  Let $\pi': \R^{n_0+r}\ra \R$ be a coordinate function and
  $u=\pi'\circ \pi_A$.
% Let $u: \R^n\ra \R$ be a linear
%function factoring through  $\pi$ ?????.
  Then  for any bounded open subset  $E$:  \newline
  (1)  $V_{Q, K}(u|_{E})=0$ for all $Q>n$  and all $K\ge K_0$;\newline
  (2)  $0<V_{n, K}(u|_{E})<\i$ for all $K\ge K_0$.

\end{Le}

\b{proof} %Let $\pi_i: \R^{n_0+k}\ra \R$ be one of the coordinate
  %functions and $u=\pi_i\circ \pi$.
 %Then for any left translation  $S_i$ of $e^{tA}S$, we have
  % $osc(u|_{S_i})=e^{\alpha t}$.
Let $P$ be a $K$-quasi-ball. Then there is a  $D_A$-ball $U$ with
$U\subset P\subset  KU$. Let $t_0=\ln(KK_0)$.   For some $t\in \R$
there is a  translate $S(t)$ of  $e^{tA}S$ and a  translate
$S(t+t_0)$  of  $e^{(t+t_0)A}S$ such that
      $\frac{1}{K_0}U\subset
S(t)\subset U$  and $KU\subset S(t+t_0)\subset KK_0 U$.
   %Then for
   Observe  that
for any  translate   $S'$ of $e^{tA}S$, we have
   $\text{osc}(u|_{S'})=e^{\la t}$.
 It follows  that
  $$\text{osc}(u|_P)\ge \text{osc}(u|_{S(t)})=\frac{1}{(KK_0)^\la}\text{osc}(u|_{S(t+t_0)})\ge
\frac{1}{(KK_0)^\la}\text{osc}(u|_P).$$  %For each quasi-ball $U$, let $V(U)$ be
  %the $n$-dimensional Lebesgue measure of the $U$.
    Also   notice that
 $\text{osc}(u|_{S(t)})=(\vol(S(t)))^{\frac{1}{n}}\le (\vol(P))^{\frac{1}{n}}$.

 Now let
  $E$ be a  bounded open subset   and
 $\{P_i\}$  a packing of  $E$ by a disjoint
 collection of $K$-quasi-balls with size $<\delta$.  For each $P_i$,
 let $U_i$ be a $D_A$-ball with $U_i\subset P_i\subset K U_i$ and
 let $S_i$ be a  translate   of some $e^{t_iA}S$ with
 $\frac{1}{K_0}U_i\subset S_i\subset U_i$. Then the preceding paragraph
      implies
 $$\sum_i \text{osc}(u|_{P_i})^Q\le (KK_0)^{Q\la} \sum_i \text{osc}(u|_{S_i})^Q\le
(KK_0)^{Q\la} \sum_i \vol(P_i)^{\frac{Q}{n}}.$$
   From this it is clear that $V_{Q, K}(u|_{E})=0$ if $Q>n$ and
$V_{n, K}(u|_{E})<\i$ since $\{P_i\}$ is a disjoint collection in
$E$.

Now consider  a particular packing $\{P_i\}$ of $E$ by the images of
the integral  unit cubes in $\R^n$ under $e^{tA}$.
 Then  $\text{osc}(u|_{P_i})=e^{\la t}$. The cardinality of  $\{P_i\}$
     is approximately
  $\frac{\vol(E)}{e^{n\la t}}$. Hence
   $V^\delta_{n, K}(u|_{E})\ge \sum_i \text{osc}(u|_{P_i})^n\approx \vol(E)$.  Hence we
   have
$0<V_{n, K}(u|_{E})<\i$.

\end{proof}

\b{Le}\label{l3.5} Suppose $A$ has only one eigenvalue $\la>0$.
  Let $E\subset \R^n$ be a rectangular box whose edges are  parallel
  to the coordinate axes.
    Let $u: (E, D_A) \ra \R$ be a
continuous function.
   Suppose there is
  some fiber $H$ of $\pi_A: \R^n\ra \R^{n_0+r}$ such that $u\vert_{H\cap E}$ is
   not constant. Then $V_{Q, K}(u)=\i$ for all $Q\le n$  and
       all $K\ge K_0$.
          %  and all sufficiently large Euclidean ball $E$. % where $K_0$ is the constant in Lemma ??.

\end{Le}

\b{proof}   Suppose   there is some
    fiber $H$ of $\pi_A$ such that $u\vert_{H\cap E}$ is
   not constant.   Then there is some Jordan block $J$ in $A$  with
   the following  property:
   if we denote  by $x=(x_1, \cdots, x_m)$ the coordinates
   corresponding to $J$, then there  is some index $k$, $1\le k\le
   m-1$ such that $u$ is
   constant along every line parallel to the $x_j$-axis   for $j\le
   k-1$,    but is
   not  constant along   some  line $L$ parallel to
   the $x_k$-axis.
  We write $\R^n=\R^{k-1}\times \R\times \R^{n-k}$,
    where the $\R$ corresponds to the $x_k$-axis and the $\R^{k-1}$
     is spanned by the $x_j$-axes ($ j\le k-1$).
        After composing $u$ with an affine function, we may assume
        that
           for some rectangular box $C=\prod_{i=1}^n [a_i, b_i]\subset E$,
              we  have
  $u\le 0$ on the codimension 1 face   $F_0:=\{x\in C: x_k=a_k\}$ of $C$
  and
   $u\ge 1$ on the codimension1 face $F_1:=\{x\in C:  x_k=b_k\}$  of $C$.
   We will induct on $k$.

 Recall that for a Jordan block $J=\la I_m +N$, we have
$e^{tJ}=e^{\la t} e^{tN}$.   %An expression for $e^{tN}$
   See (\ref{jormatrix})  for    an expression of   $e^{tN}$.

%can be found before Proposition \ref{farb-m}.

     We first assume $k=1$.
  % For any $\delta>0$, choose $t<<0$ so that the size
   %of $e^{tA}S$ is small than $\delta$.
 For $t<<0$,  consider the images of the integral unit cubes under $e^{tA}$.
   Let $\{B_i\}$ be the
   collection of all those images that intersect the box $C$.
 Notice that a vertical stack (i.e. parallel to the $x_m$-axis)  of integral cubes  is mapped by
 $e^{tA}$ to a sequences of $K_0$-quasi-balls which is almost parallel to
 the $x_1$-axis.
   We  divide   $\{B_i\}$   into
      such sequences which join
       $F_0$ and $F_1$.   Note that the projection of
         each  $B_i$
to the $x_1$-axis has length
   comparable to
 $e^{\la t}|t|^{m-1}$.  Hence  the cardinality of each sequence  is comparable to
 $e^{-\la t}|t|^{1-m}$.
    The $Q$-variation of $u$ along each  sequence is at least the
    $Q$-variation of $u$ when oscillations of $u$ on the members of
    the sequence are equal.  Since  $u\le 0$ on $F_0$ and $u\ge 1$ on
     $F_1$, this common oscillation is    at
     least  comparable to
$e^{\la t}|t|^{m-1}$.
  Since each $B_i$ has volume $e^{n\la t}$,
    the cardinality of $\{B_i\}$
    is     comparable to  $e^{-n\la t}$. It follows
  that
 the $Q$-variation of $u|_C$ is at least  comparable to
 $$e^{-n\la t}\cdot \left(e^{\la t}|t|^{m-1}\right)^Q=e^{\la t
 (Q-n)}|t|^{Q(m-1)},$$
  which $\ra \i$ as $t\ra -\i$ if $Q\le n$.
    Hence  $V_{Q, K}(u|_C)=\i$  for  $Q\le n$.

Now we assume  $m-1\ge k\ge 2$.  Then  $u$ is constant along affine
subspaces parallel to $\R^{k-1}\times \{0\}\times \{0\}\subset
\R^n$.   % and is not constant along  some  line
  %$L$ parallel to the $x_{k}$-axis.
  %We may assume $u\le 0$ on the
  %co-dimension 1 face $F_0:=\{x_{k}=0\}$ of $S$ and $u\ge 1$ on the co-dimension 1
  %face $F_1:=\{x_{k}=1\}$ of $S$.
  %Again we let ${\mathcal{B}}=\{B_i\}$ be the collection of images of
  %integral cubes under $e^{tA}$ (for $t<<0$) that intersect $S$.
   %We write $\R^n=\R^{k-1}\times \R\times \R^{n-k}$.
   %The image of the $(k-1)$-face
   %$(\R^{k-1}\times\{0\}\times \{0\})\cap S$ of $S$  under $e^{tA}$ has diameter(in Euclidean metric)
  % comparable to
    %$e^{\alpha t}|t|^{k-2}$.  The projection of $e^{tA}S$ to the
    %$x_{k}$-axis has length  comparable to $e^{\alpha t}|t|^{m-k}$.  % and the
 % projection to the subspace $\{0\}\times \{0\}\times \R^{n-k}$ has diameter
  %  comparable to
  % $e^{\alpha t}|t|^{m-k-1}$.
%Let $\text{int}(F_0)$ be the open $(n-1)$-face of $S$  with  closure
  %$F_0$.
Let $$U=\big\{x\in F_0:  (3a_i+b_i)/4\le x_i \le (a_i+3b_i)/4 \;
{\text{for all}}\;
   i\not=k
    \big\}\subset F_0.$$
 For $t<<0$, denote by $v(t)=(-1)^{m-k}e^{-\la t}e^{tA}e_m$. Notice that
 the
components  of $v(t)$  corresponding to the Jordan block $J$ is
  $$(-1)^{m-k}
 \left(\frac{t^{m-1}}{(m-1)!}, \; \frac{t^{m-2}}{(m-2)!},  \;\cdots,\; t,\; 1\right)$$
  and   all other components are $0$.  Hence  for $t<<0$,  lines
parallel to  $v(t)$ travel     much  faster  in the $x_i$ ($1\le
i\le m-1$)  direction than in the   $x_{i+1}$ direction.
 Let $Z\subset \R^n$ be the subset given by:
  $$Z=\Big\{f+ s\, v(t):  f\in U, \; 0\le s \le
  \frac{(m-k)!}{|t|^{m-k}}(b_k-a_k)\Big\}.$$
       Notice that for each fixed $f$, the
  segment $\{f+s\, v(t):0\le s \le
  \frac{(m-k)!}{|t|^{m-k}}(b_k-a_k)\}$  joins the two hyperplanes $x_k=a_k$ and
  $x_k=b_k$. Also notice that these segments are parallel to the
  images of vertical stacks (i.e., parallel to the $x_m$-axis) of integral cubes under $e^{tA}$.
   Hence $Z$ has a packing  $\mathcal{P}$
      that can be divided into sequences such
   that each sequence joins $x_k=a_k$ and $x_k=b_k$ and is the image (under $e^{tA}$)
   of a
   vertical stack of integral cubes.

For $p=(x, y, z), q=(x', y', z')\in \R^{k-1}\times \R\times
\R^{n-k}$, define $p\sim q$ if $y'=y$, $z'=z$  and $x'_i-x_i$ is an
integral multiple of $b_i-a_i$ for $1\le i\le k-1$.
% $x'-x$ have
    %integer components  and $y'=y$, $z'=z$.
     Set $Y=\R^n/\sim$ and let $\pi: \R^n\ra Y$ be the
natural
projection.  %It is easy to check that $E$ injects into $Y$.
 Also let $\pi_C: C\ra Y$ be the composition of the inclusion $C\subset \R^n$
 and $\pi$.  It is clear that $\pi_C$ is injective on the interior
 of $C$.
It is also easy to check that $\pi|_Z$ is injective.
 Now the packing $\mathcal{P}$ of $Z$  projects onto a packing of
 $Y$, which can then be pulled back through $\pi_C$ to obtain a
 packing $\mathcal{P'}$ of $C$ (since  $\pi(Z)\subset \pi_C(C)$).  A sequence in $\mathcal{P}$ gives
 rise to a broken sequence in $\mathcal{P'}$: the broken sequence will
 first hit the boundary of $C$ at a  point of $\partial
 (\prod_{i=1}^{k-1}[a_i, b_i])\times \prod_k^n [a_i, b_i]\subset \partial C$, it continues
 after a translation by an element of the form $(\sum_{i=1}^{k-1} m_i(b_i-a_i), 0,0)\in \R^n$,
 where $m_i\in \Z$;
   this can be repeated until the sequence  hits $x_k=b_k$.
 Note that we can apply Jensen's inequality to each broken sequence
 while considering $Q$-variations of $u$ since by assumption $u$ is
 constant along affine spaces parallel to $\R^{k-1}\times \{0\}\times
 \{0\}$.

%Again we organizes the members of $\mathcal B$ into sequences that
%are the images of  vertical stacks of integral cubes under $e^{tA}$.
 % Each sequence start with some $B_i$ intersecting $F_0$. But in this
  %case, the sequence will exit the unit cube $S$ at the co-dimension face
  %$x_1=0$ or $x_1=1$.  We need  sequences that exit $F_1$.  Here we
  %use the assumption that $u$ is constant along affine subspaces
  %parallel to $\R^k$.   let $B_1, \cdots, B_j$ a  the sequence
  %starting at $F_0$ and exiting at $x_1=1$. Then $B_j$ contains
  %some point $(1, a_2, \cdots, a_n)\in S$.  Then we define $B_{j+1}$
  %to be the member of $\mathcal B$ that contain the point

  Each broken sequence   joins  $F_0$ to $F_1$. Since
    the projection of $e^{tA}S$ to the
  $x_{k}$-axis has length  comparable to $e^{\la t}|t|^{m-k}$,
  the cardinality of each sequence is comparable to
    ${e^{-\la t}}{|t|^{k-m}}$.  %Similar to the arguments in the  third  paragraph (of the proof),
   The $Q$-variation of $u$
  along the sequence is at least the $Q$-variation when the
  oscillations of $u$ are the same on all members of the sequence.
   The common oscillation is at least comparable to $e^{\la t}|t|^{m-k}$.
  Hence the $Q$-variation of $u$ is at  least  comparable to
    $$\frac{1}{e^{n\la t}}   \cdot  (e^{\la
    t}|t|^{m-k})^Q=|t|^{Q(m-k)}e^{(Q-n)\la t},$$
 which $\ra \i$ when $t\ra -\i$ if $Q\le n$. Hence
  $V_{Q, K}(u|_C)=\i$ for $Q\le n$.

\end{proof}

\section{Proof of the main  theorems}\label{pro}

In this section we  prove the main results of the paper.
 The main tools are the notion of $Q$-variation (Section \ref{qvar})
    and the arguments
 from Section 4 of \cite{X}  and \cite{SX}.  The main results of \cite{SX} and
 \cite{X} are the basic steps in the induction.

We first fix the notation. Let $A$ be an $n\times n$ matrix in real
part Jordan form with positive eigenvalues
 $$\la_1<\cdots<\la_{k_A}.$$
  Let $V_i\subset\R^n$ be the generalized eigenspace of $\la_i$, and
set $d_i=\dim V_i$.
   If $k_A\ge 2$, we write $A$ in the block diagonal form
    $A=[A_1, \cdots, A_{k_A}]$,  where $A_i$ is the block
    corresponding to the eigenvalue  $\la_i$;
 we also denote $A'=[A_1, \cdots, A_{k_A-1}]$.  % and $B'=[B_1, \cdots, B_{k_B-1}]$.
   If  $k_A=1$, that is, if $A$ has only one eigenvalue
   $\la=\la_1$, then we also write  $A=[\la I_{n_0}, \la I_{n_1}+N, \cdots,
   \la I_{n_r}+N]$ in the block diagonal form, and we let
    $\pi_A: \R^n\ra \R^{n_0+r}$ be the projection
       %of
    %$\R^n$ onto the eigenspace of $A$
        defined before   Lemma \ref{l11}. If $k_A=1$ and $r\ge 1$, we set
    $l_A=\max\{n_1,  \cdots, n_r\}$.

  Similarly,  let  $B$ be an $n\times n$ matrix in
real part Jordan form with positive eigenvalues
 $$\mu_1<\cdots<\mu_{k_B}.$$
Let $W_j\subset\R^n$ be the generalized eigenspace of $\mu_j$, and
set $e_j=\dim W_j$.
   If $k_B\ge 2$, we write $B$ in the block diagonal form
    $B=[B_1, \cdots, B_{k_B}]$,  where $B_j$ is the block
    corresponding to the eigenvalue  $\mu_j$;
     we also denote  $B'=[B_1, \cdots, B_{k_B-1}]$.
  If $k_B=1$, that is, if $B$ has only one eigenvalue
   $\mu=\mu_1$,  we also write  $B=[\mu I_{m_0}, \mu I_{m_1}+N, \cdots,
   \mu I_{m_s}+N]$ in the  block diagonal form, and we let
    $\pi_B: \R^n\ra \R^{m_0+s}$ be the   projection
        %of
      %$\R^n$ onto the eigenspace of $B$
         defined before   Lemma \ref{l11}.  If $k_B=1$  and $s\ge 1$, we set
     $l_B=\max\{m_1, \cdots, m_s\}$.

     Suppose there is an   $\eta$-quasisymmetric map
 $F: (\R^n, D_A)\ra (\R^n, D_B)$.

\b{Le}\label{l4.1} $k_A=1$ if and only if $k_B=1$.

\end{Le}

\b{proof} Suppose $k_A=1$ and $k_B\ge 2$.
   Fix any $Q$ with
 $\frac{\sum_j{\mu_j}e_j}{\mu_{k_B}}  < Q< n$.
Let $\pi: (\R^n, D_B)\ra W_{k_B}$ be the projection onto $W_{k_B}$,
and
  $\pi': W_{k_B}\ra \R$  a coordinate function on $W_{k_B}$.  Set
  $u=\pi'\circ \pi$.
  %Let $u: \R^{n}\ra
 %\R$ be any nontrivial linear function  factoring through
 % the projection to $V_{r'}$.
    Then Lemma \ref{l3.3} implies $V_{Q,
  \eta(K)}(u|_{F(E)})=0$ for all  sufficiently   large  $K$  and all   bounded subsets  $E\subset(\R^n, D_A)$.  By Lemma \ref{l3.0} $V_{Q, K}(u\circ
  F|_E)=0$.  But this contradicts Corollary \ref{c3.1}.

\end{proof}

\b{Le}\label{l4.2}
  Suppose $k_A=1$.  Then $A=\la_1 I_n$ if and only if $B= \mu_1
  I_n$.

\end{Le}

\b{proof} Suppose   $B= \mu_1
  I_n$.
     Let $\pi_i$ ($i=1, 2, \cdots, n$) be the coordinate functions on
     $(\R^{n}, D_B)$. Then  by Lemma \ref{l3.4}  we have $V_{n,
     \eta(K)}(\pi_i|_{F(E)})<\i$ for all  $i$,  all sufficiently large  $K$
       and all rectangular boxes  $E\subset (\R^n, D_A)$. Hence
      $V_{n, K}(\pi_i\circ F|_{E})<\i$ by Lemma \ref{l3.0}.
        % Let $\pi$ be the canonical projection of $\R^n$ onto the
         %eigenspace of $A$, as defined before Lemma \ref{l3.4}.
       Now Lemma \ref{l3.5}  implies  that $\pi_i\circ F$
        is constant on the fibers of $\pi_A$. Since this is true for
        all $1\le i\le n$,  the fibers of $\pi_A$ must have  dimension
        $0$. Hence $A$ must also be a multiple of $I_n$.

\end{proof}

\b{Le}\label{l4.3} Suppose $k_A=1$ and $r\ge 1$. Then
 $F$ maps each fiber of $\pi_A$ onto some fiber of $\pi_B$.

\end{Le}

\b{proof} Lemmas \ref{l4.1} and  \ref{l4.2} imply that  $k_B=1$ and
$s\ge 1$.
  Notice that it suffices to show that each fiber of $\pi_A$ is
  mapped by $F$ into some fiber of $\pi_B$:  by symmetry each fiber
  of $\pi_B$ is mapped by $F^{-1}$ into some fiber of $\pi_A$ and
 hence  the lemma follows.
We shall prove this by contradiction and so  assume that there is
some fiber $H$ of $\pi_A$ such that $F(H)$ is not contained in any
fiber of $\pi_B$.  Then there is some  coordinate  function
 $\pi': \R^{m_0+s}\ra \R$  such that  $u\circ F$ is
    not constant on $H$,
      where $u:=\pi'\circ \pi_B$. Now Lemma \ref{l3.4} implies
   that  $V_{n, \eta(K)}(u|_{F(E)})<\i$  for all sufficiently  large $K$
     and all  rectangular boxes  $E\subset (\R^n, D_A)$. By Lemma
   \ref{l3.0} we have $V_{n, K}(u\circ F|_{E})<\i$.
This contradicts Lemma \ref{l3.5} since we can choose $E$ such that
  $u\circ F$ is not constant on $H\cap E$.

\end{proof}

It follows from Lemma \ref{l4.3} that
 $F$ induces a map
 $G:  \R^{n_0+r}\ra \R^{m_0+s}
    $
    such that
  $F(\pi_A^{-1}(y))=\pi_B^{-1}(G(y))$
     for all   $y\in \R^{n_0+r}$.
  Define
$$\tau_A: \R^n=\R^{n_0}\times \R^{n_1}\times \cdots\times \R^{n_r}\longrightarrow \R^n=\R^{n-n_0-r}\times \R^{n_0+r}$$
   by
      $$\tau_A(z, (x_1, y_1), \cdots, (x_r, y_r))= ((x_1, \cdots, x_r),   (z, y_1, \cdots,
      y_r)),$$
          where $(x_i, y_i)\in
      \R^{n_i}=\R^{n_i-1}\times \R$.  Similarly, there is an
      identification   $$\tau_B:  \R^n=\R^{m_0}\times \R^{m_1}\times \cdots\times \R^{m_s}\longrightarrow \R^n=\R^{n-m_0-s}\times
      \R^{m_0+s}.$$
  With the identifications $\tau_A$  and $\tau_B$, we have
$\pi_A^{-1}(y)= \R^{n-n_0-r}\times \{y\}$,  $\pi_B^{-1}(G(y))=
\R^{n-m_0-s}\times \{G(y)\}$,
    and $F(\R^{n-n_0-r}\times \{y\})=\R^{n-m_0-s}\times \{G(y)\}$.
%      Now $\pi_A^{-1}(y)$ is identified with $\{y\}\times \R^{n-n_0-r}$
 % via $\tau$.
 % We identify $\pi_A^{-1}(y)$ with $\R^{n-n_0-r}$ and
   %$\pi_B^{-1}(G(y))$  with  $\R^{n-m_0-s}$.
   Hence for each $y\in
   \R^{n_0+r}$, there is a  map
   $$H(\cdot, y):\R^{n-n_0-r}\ra\R^{n-m_0-s}$$
     such that $F(x,y)=(H(x,y), G(y))$  for all $x\in\R^{n-n_0-r}$.

\b{Le}\label{l4.35} {Suppose $k_A=1$ and $r\ge 1$.
  Then:\newline
  (1)  The map $G:  (\R^{n_0+r}, |\cdot|^{\frac{1}{\la}})\ra (\R^{m_0+s},
    |\cdot|^{\frac{1}{\mu}})$ is  $\eta$-quasisymmetric;\newline
 %    $G:  (\R^{n_0+r}, |\cdot|^{\frac{1}{\la}})\ra (\R^{m_0+s},
  %  |\cdot|^{\frac{1}{\mu}})$ such that
     %$F(\pi_A^{-1}(y))=\pi_B^{-1}(G(y))$;\newline
 (2)  for each $y\in  \R^{n_0+r}$, the   map $H(\cdot, y):  (\R^{n-n_0-r}, D_{A(1)})\ra
    (\R^{n-m_0-s}, D_{B(1)})$
  is $\eta$-quasisymmetric.
    %$F:  (\pi_A^{-1}(y),  D_{A(1)})\ra (\pi_B^{-1}(G(y)),
    % D_{B(1)})$ is $\eta$-quasisymmetric, where
     %  $\pi_A^{-1}(y)$
    %and $\pi_B^{-1}(G(y)$  are
     %   naturally identified with $\R^{n-n_0-r}$
     %  and  $\R^{n-m_0-s}$ respectively.

}

\end{Le}

\b{proof} (1) follows from Lemma \ref{l12}  and the arguments on
page 10 of \cite{X}.   The statement (2) follows from Lemma
\ref{l11}.

\end{proof}

 Suppose $k_A=1$.  Set  $\epsilon=\la/\mu$ and
   $\eta_1(t)=\eta(t^{\frac{1}{\epsilon}})$. We notice that all the
   following maps are $\eta_1$-quasisymmetric:\newline
   (1) $F: (\R^n, D^\epsilon_A)\ra (\R^n, D_B)$;\newline
   (2)  $G:  (\R^{n_0+r}, |\cdot|^{\frac{1}{\mu}})\ra (\R^{m_0+s},
    |\cdot|^{\frac{1}{\mu}})$;\newline
    (3)  $H(\cdot, y):  (\R^{n-n_0-r}, D^\epsilon_{A(1)})\ra
    (\R^{n-m_0-s}, D_{B(1)})$, for each
$y\in  \R^{n_0+r}$.

Let $g: (X_1, \rho_1)\ra (X_2, \rho_2)$ be a  bijection  between two
quasimetric spaces.
   Suppose $g$ satisfies the following condition:
      for any  fixed $x\in X_1$,  $\rho_1(y,x)\ra  0$
           if and only if
      $\rho_2(g(y), g(x))\ra 0$.
We define for every $x\in X_1$  and $r>0$,
\begin{align*}
   L_g(x,r)&=\sup\{\rho_2(g(x), g(x')):   \rho_1(x,x')\le r\},\\
   l_g(x,r)&=\inf\{\rho_2(g(x), g(x')):   \rho_1(x,x')\ge r\},
\end{align*}
and set
\[
   L_g(x)=\limsup_{r\ra 0}\frac{L_g(x,r)}{r}, \ \
   l_g(x)=\liminf_{r\ra 0}\frac{l_g(x,r)}{r}.
\]

\begin{Le}\label{l36}   Consider  the maps
$G:  (\R^{n_0+r}, |\cdot|^{\frac{1}{\mu}})\ra (\R^{m_0+s},
    |\cdot|^{\frac{1}{\mu}})$   and
      $H(\cdot, y):  (\R^{n-n_0-r}, D^\epsilon_{A(1)})\ra
    (\R^{n-m_0-s}, D_{B(1)})$.
  The following hold for all $y\in \R^{n_0+r}$, $x\in \R^{n-n_0-r}$:
\newline
(1)  $L_G(y, r)\le\eta_1(1)\, l_{H(\cdot, y)}(x, r)$ for any
$r>0$;\newline
   (2)  $\eta_1^{-1}(1)\, l_{H(\cdot, y)}(x)\le l_G(y)\le \eta_1(1)\, l_{H(\cdot, y)}(x)$; \newline
    (3)  $\eta_1^{-1}(1)\, L_{H(\cdot, y)}(x)\le L_G(y)\le \eta_1(1)\, L_{H(\cdot,
    y)}(x)$.

   % $L_G(y)\ge  \eta^{-1}(1) L_{H(\cdot, y)}(x)$.

   \end{Le}

\b{proof}
  The proof is very similar to that of Lemma 4.3 in \cite{X}.
  Let $y\in \R^{n_0+r}$, $x\in \R^{n-n_0-r}$  and  $r>0$.
 Let $y'\in \R^{n_0+r}$   with
 $|y-y'|^{\frac{1}{\mu}}\le r$  and     $x'\in \R^{n-n_0-r}$
   with $D^\epsilon_{A(1)}(x, x') \ge r$.   Set $t_0=\ln |y'-y|/\la$.
     Let $(u_i, v_i)$ ($u_i\in \R^{n_i-1}$, $v_i\in \R$,  $1\le i\le r$) be the unique solution of
      $e^{-t_0 N_{n_i}}(u_i, v_i)^T=(0,\cdots, 0, y'_i-y_i)^T$.
        Let $x''_i=u_i+x_i$ and $x''=(x''_1, \cdots, x''_r)$.
  Then $D^\epsilon_A((x,y), (x'', y'))=|y-y'|^{\frac{1}{\mu}} \le r\le D^\epsilon_A((x, y), (x', y))$.
    Since $F: (\R^n, D^\epsilon_A)\ra (\R^n, D_B)$ is $\eta_1$-quasisymmetric, we have
\b{align*} |G(y)- G(y')|^{\frac{1}{\mu}}   \le D_B(F(x'',y'),
F(x,y)) & \le \eta_1(1)\,
D_B(F(x,y), F(x', y))\\
& =\eta_1(1)\,D_{B(1)}(H(x,y),H(x', y)).
\end{align*}
    Since $y'$ and $x'$ are chosen arbitrarily,    (1)  follows.

   The proofs of (2) and (3) are exactly the same as those for Lemma
   4.3 in \cite{X}.

\end{proof}

% If $A=I_2+N$ is the  $2\times 2$ Jordan block with eigenvalue 1,
% then $D_A$ is biLipschitz to the following function  (see
% Proposition 2.2 in \cite{X}):
%  \begin{equation}\label{2times2}
 % D((x_1, y_1)^T, (x_2, y_2)^T)=\max \{|y_2-y_1|,
  %|(x_2-x_1)-(y_2-y_1)\ln |y_2-y_1||\},
  %\end{equation}
  %where $0\ln 0$ is understood to be $0$.

\b{Le}\label{l4.4} Suppose $k_A=1$ and $l_A=2$. Then $l_B=2$ and
   for $\epsilon=\la/\mu$\newline
  %there is  a  constant $\epsilon>0$ such that:\newline
  (1) $A$ and
$\epsilon B$ have the same real part Jordan form;\newline
  (2)  The   map
  $F: (\R^n, D_A^\epsilon)\ra (\R^{n}, D_B)$  is a  $K$-quasisimilarity, where
  $K$ depends only on  $A$, $B$ and  $\eta$.  %  and  $\epsilon$.

\end{Le}

\b{proof}  %Assume $l_A=2$.  %We write $p, p'\in (\R^n, D_A)$ as
%$p=(z, (x_1, y_1), \cdots, (x_r, y_r))^T$
 % and   $p'=(z',
 %(x'_1, y'_1), \cdots, (x'_r, y'_r))^T$
  % as indicated before Lemma
%\ref{l3.4}.  Then $D_A(p, p')$ is biLipschitz to
 %$$\max\{|z'-z|, |y'_i-y_i|, |(x'_i-x_i)-(y'_i-y_i)\ln
  %|y'_i-y_i||\}.$$
 % Notice that $\pi_A(p)=(z, y_1, \cdots, y_r)$.  It follows that
    %This together with  the formula  implies that
  %each fiber of $\pi_A$
  %is biLipschitz to a Euclidean space.
 % Let $H$ be a fiber of $\pi_A$.  By Lemma \ref{l4.3} $F(H)$ is a
  %fiber of $\pi_B$.   Notice that
   % formula (\ref{jormatrix}) implies that the restriction of $D_A$ on
   % $H$ agrees with $D_{A(1)}$, where $A(1)=[\la I_{n_1-1}+N,
    %\cdots, \la I_{n_r-1}+N]$ and $\la I_{n_i-1}+N$ is understood to be
    %$\la I_1$ if $n_i=2$.  Similarly for $F(H)$.
   (1)   By Lemma \ref{l4.35} (2),
 for each $y\in  \R^{n_0+r}$, the   map $H(\cdot, y):  (\R^{n-n_0-r}, D_{A(1)})\ra
    (\R^{n-m_0-s}, D_{B(1)})$
  is $\eta$-quasisymmetric.
    Since $l_A=2$,  all Jordan blocks of  $A$ have size $2$ and  $A(1)=\la I_r$.  % is a  multiple of the identity matrix
            %$I_r$.
    Now  Lemma \ref{l4.2} applied to  $H(\cdot, y)$
  implies that
  $B(1)=\mu I_r$.    % is also a multiple of the identity matrix  $I_r$.
%  each fiber of $\pi_B$ is also
%biLipschitz to a Euclidean space.
  It follows that all Jordan blocks of $B$ also have size $2$, and
    hence $l_B=2$ and $B(1)=\mu I_s$.  So we have $r=s$. That is, $A$ and $B$ have the same number of $2\times 2$ Jordan blocks.
     Now  (1) follows.

%  It then follows that $B(1)$

 % Lemma \ref{l4.3} also implies that the fibers of $\pi_A$ and
  %$\pi_B$ have the same  dimension.  It follows that $A$ and $B$
  %have the same number of $2\times 2$ Jordan blocks.  %Since we also
   % Now  (1)   holds with $\epsilon=\la/\mu$.
(2)   The proof of (2) is very similar to the arguments in Section 4
of \cite{SX}  and
  \cite{X}.   We will only indicate the differences here.
    First we notice that  $G:
(\R^{n_0+r}, |\cdot|)\ra (\R^{m_0+s},
    |\cdot|)$  is also quasisymmetric, and hence is
    differentiable a.e.
    Since  $F: (\R^n, D^\epsilon_A)\ra (\R^n, D_B)$ is
   $\eta_1$-quasisymmetric,
   the
  arguments in Section 4 of
  \cite{SX}  and \cite{X}  imply that there is a constant $K_1$ depending only on
$\eta_1$, such that for every $y\in \R^{n_0+r}$ where $G:
(\R^{n_0+r}, |\cdot|)\ra (\R^{m_0+s},
    |\cdot|)$
  is differentiable,   we have $0<l_G(y)<\i$  and  the   map
  $$H(\cdot, y):  (\R^{n-n_0-r}, D^\epsilon_{A(1)})\ra
    (\R^{n-m_0-s}, D_{B(1)})$$
        is  a  $K_1$-quasisimilarity with
    constant $l_G(y)$.

    Now  let $y, y'\in \R^{n_0+r}$
  be two  points where $G$ is differentiable. We will show that
  $l_G(y)$ and $l_G(y')$ are comparable.
Let $x \in \R^{n-n_0-r}$ and choose
   $x'\in \R^{n-n_0-r}$ so that $D_{A(1)}(x, x')>>|y'-y|^{\frac{1}{\la}}$.
Let $(u_i, v_i)$ be as in the proof
  of Lemma  \ref{l36}. Let $x''_i=x_i+u_i$, $x'''_i=x'_i+u_i$ ($1\le i\le
  r$), and set  $x''=(x''_1, \cdots, x''_r)$, $x'''=(x'''_1, \cdots,
  x'''_r)$.  Then $$D_A((x,y), (x'', y'))=D_A((x', y), (x''',
  y'))=|y'-y|^{\frac{1}{\la}}.$$
        Now   the generalized triangle inequality implies
\b{align*}D_A((x'',y'), (x', y)) &\le M\Big\{D_A((x'',y'),
(x,y))+D_A((x,y),
       (x', y))\Big\}\\
       & \le 2MD_A((x,y),
       (x', y)).
       \end{align*}
       By the quasisymmetry condition we have
       $$D_B(F(x'',y'),
       F(x',y))\le \eta(2M) D_B(F(x,y), F(x', y)).$$
         Similarly, $D_B(F(x'',y'), F(x''', y'))\le \eta(2M)D_B(F(x'',y'),
       F(x',y))$. So we have
$$D_B(F(x'',y'), F(x''', y'))\le (\eta(2M))^2 D_B(F(x,y), F(x',
y)).$$
   This together with the quasisimilarity properties of
$H(\cdot,y)$ and $H(\cdot, y')$ mentioned above implies   that
$$l_G(y')D^\epsilon_{A(1)}(x'', x''')\le K_1^2(\eta(2M))^2 l_G(y)
D^\epsilon_{A(1)}(x, x').$$
     Since $D_{A(1)}(x'', x''')=D_{A(1)}(x,
x')$, we have $l_G(y')\le K_1^2(\eta(2M))^2 l_G(y)$.  By symmetry,
we also have $l_G(y)\le K_1^2(\eta(2M))^2 l_G(y')$.
  Now fix $y$ and set $C=l_G(y)$.
 Then  at every $y'$ where $G$ is differentiable,
  $H(\cdot,y')$ is  a $K_2$-quasisimilarity with constant $C$,
  where $K_2= K_1^3(\eta(2M))^2$.  Now   a limiting argument shows
  that this is true for every $y'\in \R^{n_0+r}$.
 The arguments in
Section 4 of
  \cite{X}
   (using Lemma \ref{l36} from above instead of Lemma  4.3  in
   \cite{X})
   then show that there is a constant $K_3=K_3(K_2, \eta_1)$ such that
    $G:  (\R^{n_0+r}, |\cdot|^{\frac{1}{\mu}})\ra (\R^{m_0+s},
    |\cdot|^{\frac{1}{\mu}})$
  and all $H(\cdot, y)$ are $K_3$-quasisimilarities    with constant
  $C$.

   The final difference is in finding a lower bound for $D_B(F(x,y),
   F(x',y'))$.
  %Let $(u_i, v_i)$ be as in the above paragraph.
   If $D^\epsilon_A((x,y), (x',y'))\le (2M)^\epsilon |y'-y|^{\frac{1}{\mu}}$, then
\b{align*}
D_B(F(x,y),
   F(x',y'))\ge |G(y')-G(y)|^{\frac{1}{\mu}}& \ge \frac{C}{K_3}
   |y'-y|^{\frac{1}{\mu}}\\
   & \ge
\frac{C}{(2M)^\epsilon K_3}
   D^\epsilon_A((x,y), (x',y')).
   \end{align*}
   Now assume
$D^\epsilon_A((x,y), (x',y'))\ge (2M)^\epsilon
|y'-y|^{\frac{1}{\mu}}$. Let $(u_i, v_i)$ be as in the above
paragraph. Let $x''_i=x'_i-u_i$ and set $x''=(x''_1, \cdots,
x''_r)$. Then
 $D^\epsilon_A((x'', y), (x',y'))=|y'-y|^{\frac{1}{\mu}}$.
    The generalized triangle inequality implies
      $$\frac{1}{2M}\le \frac{D_A((x,y), (x'', y))}{D_A((x,y), (x',
      y'))}\le 2M.$$
      Now the quasisymmetric condition implies
\b{align*} D_B(F(x,y), F(x',y'))& \ge \frac{1}{\eta(2M)}D_B(F(x,y),
F(x'', y))\\
&  \ge \frac{C}{K_3\eta(2M)} D_A^\epsilon((x,y), (x'',
y))\\
& \ge \frac{C}{(2M)^\epsilon K_3\eta(2M)}D^\epsilon_A((x,y), (x',
      y')).
      \end{align*}
        So we have found a lower bound for
 $D_B(F(x,y),
   F(x',y'))$.  The rest of the proof is the same as in
Section 4 of
  \cite{X}.  We notice that the constant $M$ depends only on $A$,
   and $\epsilon$ depends only on $A$ and $B$.
    Hence
$F$ is a $K$-quasisimilarity with $K$ depending only on $A$, $B$ and
$\eta$.

\end{proof}

\b{Le}\label{l4.5} Suppose $k_A=1$ and $l_A\ge 2$. Then for
$\epsilon=\la/\mu$:\newline
      (1)   $A$ and $\epsilon B$ have the same
real part Jordan form;\newline
  (2)  The   map
  $F: (\R^n, D^\epsilon_A)\ra (\R^{n}, D_B)$  is a  $K$-quasisimilarity, where
  $K$ depends only on   $A$, $B$  and $\eta$.

\end{Le}

\b{proof}
 We induct on $l_A$. The basic step  $l_A=2$ is Lemma \ref{l4.4}.
  Now assume $l_A=l\ge 3$ and that the lemma holds for $l_A=l-1$.
    For any   $y\in  \R^{n_0+r}$,   the induction hypothesis
    applied to
the
   $\eta$-quasisymmetric
 map $H(\cdot, y):  (\R^{n-n_0-r}, D_{A(1)})\ra
    (\R^{n-m_0-s}, D_{B(1)})$  implies
  %  is $\eta$-quasisymmetric.
  % By Lemma \ref{l4.3}, each fiber $H$ of $\pi_A$ is mapped by $F$ onto a
 %fiber $H'$ of $\pi_B$.  The induction hypothesis applied to $F|_H:
 %H\ra H'$ implies
 that  for $\epsilon=\la/\mu$:  \newline
   %there is a  constant $\epsilon>0$ such that:\newline
  (a)   $A(1)$ and $\epsilon B(1)$
have the same real part Jordan form;\newline
 (b)  $H(\cdot, y):  (\R^{n-n_0-r}, D^\epsilon_{A(1)})\ra
    (\R^{n-m_0-s}, D_{B(1)})$  is a $K$-quasisimilarity with $K$ depending  only on
   $A(1)$, $B(1)$ and $\eta$.\newline
   Now (1) follows from (a), and  (2)
follows from (b), Lemma \ref{l36}
      and the
  arguments in Section  4 of \cite{X}
    (see the proof of Lemma \ref{l4.4}(2)).

%   follows from (b), Lemma \ref{l36}
 %      and the
  % arguments in Section  4 of \cite{X}.

\end{proof}

\b{Le}\label{l4.6} Suppose $k_A\ge 2$.
  Then $k_B\ge 2$
    and
         $\frac{\sum_i\,d_i\la_i}{\la_{k_A}}=\frac{\sum_j\,e_j\mu_j}{\mu_{k_B}}$.

\end{Le}

\b{proof} Lemma \ref{l4.1} implies $k_B\ge 2$.
  Suppose  $\frac{\sum_i\,d_i\la_i}{\la_{k_A}}>\frac{\sum_j\,e_j\mu_j}{\mu_{k_B}}$.
    Pick any $Q$ with
$\frac{\sum_i\,d_i\la_i}{\la_{k_A}}>Q>\frac{\sum_j\,e_j\mu_j}{\mu_{k_B}}$.
Let $\pi: (\R^n, D_B)\ra W_{k_B}$ be the projection onto $W_{k_B}$,
and
  $\pi': W_{k_B}\ra \R$  a coordinate function on $W_{k_B}$.  Set
  $u=\pi'\circ \pi$.
    By Lemma \ref{l3.3} we have $V_{Q, \eta(K)}(u|_{F(E)})=0$ for
  all  sufficiently large $K$  and all Euclidean balls $E\subset (\R^n, D_A)$.  Lemma \ref{l3.0} implies $V_{Q, K}(u\circ F|_E)=0$.
  This contradicts Lemma \ref{l3.1}   since $Q<\frac{\sum_i\,d_i\la_i}{\la_{k_A}}$  and  the function $u\circ F$ is
  nonconstant.  Similarly there is a contradiction if
     $\frac{\sum_i\,d_i\la_i}{\la_{k_A}}<\frac{\sum_j\,e_j\mu_j}{\mu_{k_B}}$.
    The lemma follows.

\end{proof}

Recall that (see Section \ref{metriconb}),  if $k_A\ge 2$, then the
restriction of $D_A$ to each affine subspace $H$ parallel to
$\prod_{i<k_A} V_i$  agrees with $D_{A'}$, where $A'=[A_1, \cdots,
A_{k_A-1}]$.

\b{Le}\label{l4.7} Denote $k=k_A$ and $k'=k_B$. Suppose $k\ge 2$.
Then each affine subspace $H$
 of    $\R^n$   parallel to $\prod_{i<k} V_i$ is mapped by  $F$
   onto an affine
subspace parallel to $\prod_{j<k'}W_j$. Furthermore,  $F|_H: (H,
D_{A'}) \ra (F(H), D_{B'})$ is $\eta$-quasisymmetric,  and
   $F$ induces an $\eta$-quasisymmetric map $G: (V_{k}, D_{A_{k}})\ra (W_{k'}, D_{B_{k'}})$
   such that $F((\prod_{i<k} V_i)\times \{y\})=(\prod_{j<k'}
   W_j)\times \{G(y)\}$.

\end{Le}

\b{proof}
  As in the proof of Lemma \ref{l4.3}, to establish the first claim it suffices to show that each affine
  subspace parallel to  $\prod_{i<k} V_i$ is mapped into an affine
subspace parallel to $\prod_{j<k'}W_j$. By Lemma \ref{l4.6} we have
$\frac{\sum_i\,d_i\la_i}{\la_{k}}=\frac{\sum_j\,e_j\mu_j}{\mu_{k'}}$.
  Pick any $Q$ with
    $$\frac{\sum_i\,d_i\la_i}{\la_{k}}<Q<\min\bigg\{\frac{\sum_i\,d_i\la_i}{\la_{k-1}},
    \frac{\sum_j\,e_j\mu_j}{\mu_{k'-1}}\bigg\}.$$
   Suppose there is an  affine
  subspace  $H$ parallel to  $\prod_{i<k} V_i$  such that $F(H)$ is not contained in any   affine
subspace parallel to $\prod_{j<k'}W_j$.
  Let  $\pi: \prod_j W_j\ra W_{k'}$ be  the canonical projection.
  Then there is some coordinate function $\pi': W_{k'}\ra \R$ such
  that
%Then there is a linear function $u:(\R^{n}, D_B)\ra \R$  factoring
 %through the projection onto $V_{k_b}$
     $u\circ F$ is not constant on $H$,   where $u=\pi'\circ
     \pi$.  % and $\pi: \sum_j W_j\ra W_{k_B}$ is the canonical projection.
    As  $Q>\frac{\sum_j\,e_j\mu_j}{\mu_{k'}}$,
   Lemma \ref{l3.3} implies  $V_{Q, \eta(K)}(u|_{F(E)})=0$ for all  sufficiently large  $K$
     and all rectangular boxes   $E\subset (\R^n, D_A)$.   By  Lemma \ref{l3.0}   $V_{Q, K}(u\circ F|_E)=0$.  This
  contradicts  Lemma \ref{l3.2} since
  $Q<\frac{\sum_i\,d_i\la_i}{\la_{k-1}}$  and  we can choose
  a  rectangular  box $E$ such that
   $u\circ F$ is not
  constant on $H\cap E$.

 Since  by  assumption  $F$ is $\eta$-quasisymmetric,
  it follows from the remark preceding the lemma that
  %  and the assumption that $F$ is $\eta$-quasisymmetric  that
 $F|_H: (H,
D_{A'}) \ra (F(H), D_{B'})$ is $\eta$-quasisymmetric.

 % Since $F$ is $\eta$-quasisymmetric, its restriction $F|_H: H\ra
  %F(H)$ is also $\eta$-quasisymmetric.

The first claim implies that there is a map $G:V_{k}\ra W_{k'}$
   such that $F((\prod_{i<k} V_i)\times \{y\})=(\prod_{j<k'}
   W_j)\times \{G(y)\}$   for any $y\in V_k$.  That
   $G: (V_{k}, D_{A_{k}})\ra (W_{k'}, D_{B_{k'}})$ is $\eta$-quasisymmetric follows
    from
     (\ref{eq:1}), (\ref{eq:2})  %the fact that $D_A$ and $\tilde D$ are biLipschitz equivalent  and
    and  the arguments on page 10  of \cite{X}.

\end{proof}

\b{Le}\label{l4.8} Suppose $k_A=2$.  Then $k_B=2$  and for
$\epsilon=\la_1/\mu_1$:\newline
  %   there is a  constant $\epsilon>0$ such that:\newline
       (1) $A$ and $\epsilon B$ have the same real part Jordan form;\newline
  (2)  The   map
  $F: (\R^n, D^\epsilon_A)\ra (\R^{n}, D_B)$  is a   $K$-quasisimilarity, where
  $K$ depends only on  $A$,  $B$ and $\eta$.

\end{Le}

\b{proof} Let  $H$   be an   affine subspace   of $\R^n$ parallel to
$\prod_{i<k_A} V_i$.
   By Lemma \ref{l4.7}
   $F(H)$  is  an affine subspace parallel to
$\prod_{j<k_B}W_j$,  and   $F|_H: (H, D_{A'}) \ra (F(H), D_{B'})$ is
$\eta$-quasisymmetric.
 Since $k_A=2$, we have $k_{A'}=1$.
 %Set $A'=[A_1, \cdots, A_{k_A-1}]$ and $B'=[B_1,
   %\cdots, B_{k_B-1}]$.
%Notice that $(H, D_A)$ is isometric to
 % $(\sum_{i<k_A} V_i, D_{A'})$ and
  %  similarly  $(F(H), D_B)$ is isometric to
  %$(\sum_{j<k_B} W_j, D_{B'})$.
  Now   Lemma \ref{l4.1}    applied to
   $F|_H$ implies $k_{B'}=1$, so $k_B=k_{B'}+1=2$.
 Now the $\eta$-quasisymmetric map
    $F|_H: (H, D_{A'}) \ra (F(H), D_{B'})$
  becomes $(V_1, D_{A_1})\ra (W_1, D_{B_1})$, and Lemmas \ref{l4.5}
and \ref{l4.2}  imply that $A_1$ and $\epsilon_1 B_1$ have the same
real part Jordan form, where $\epsilon_1=\la_1/\mu_1$.
   By Lemma \ref{l4.7} $F$  induces an $\eta$-quasisymmetric map
    $G: (V_2,  D_{A_2}) \ra (W_2, D_{B_2})$,
   and hence  Lemmas \ref{l4.5}  and \ref{l4.2}   again imply that $A_2$ and $\epsilon_2 B_2$ have the
same real part Jordan form, where $\epsilon_2=\la_2/\mu_2$. Lemma
\ref{l4.7}  also implies   $d_1=e_1$  and   $d_2=e_2$.  Now  Lemma
\ref{l4.6}
  implies $\la_1/\mu_1=\la_2/\mu_2$.  Hence (1) holds.

 To prove (2), we  consider two cases. First assume that %$A$ has the
   %form  $A=[\alpha_1 I_{d_1}, \alpha_2 I_{d_2}]$.
    $A_1=\la_1 I$ and $A_2=\la_2 I$.
 In this case, (2) follows from (1) and  Theorem \ref{sx}.
  Next we assume that either  $A_1\not=\la_1 I$ or $A_2\not=\la_2 I$  holds.
  Then Lemma \ref{l4.5} implies that either $F|_H: (H, D^\epsilon_{A_1})\ra (F(H), D_{B_1})$
  is a  $K_1$-quasisimilarity with $K_1$ depending  only on  $A_1$, $B_1$ and
  $\eta$,
  or   $G: (V_2, D^\epsilon_{A_2}) \ra (W_2, D_{B_2})$ is  a
  $K_2$-quasisimilarity  with   $K_2$ depending  only on
  $A_2$, $B_2$ and  $\eta$.    %  and $\epsilon$.
  Then the arguments  similar to those in the proof of Lemma \ref{l4.4}(2) (also compare with Section 4
  of \cite{SX})   show that
$F: (\R^n, D^\epsilon_A)\ra (\R^{n}, D_B)$
   is a  $K$-quasisimilarity with $K$ depending only on
    $A$, $B$ and $\eta$.

\end{proof}

\b{Le}\label{l4.9} Suppose $k_A\ge 2$.  Then
  for
$\epsilon=\la_1/\mu_1$:\newline
   (1)  $A$ and $\epsilon B$ have the same real part Jordan form;\newline
  (2)  The   map
  $F: (\R^n, D^\epsilon_A)\ra (\R^{n}, D_B)$  is a $K$-quasisimilarity, where
  $K$ depends only on  $A$, $B$ and $\eta$.

\end{Le}

\b{proof}
 We induct on $k_A$.  The basic step $k_A=2$ is Lemma \ref{l4.8}.
   Now we assume $k_A=k\ge 3$ and that the lemma holds for
    $k_A=k-1$.
 For each affine subspace $H$
of $\R^n$ parallel to $\prod_{i<k_A} V_i$, the induction hypothesis
applied to   $F|_H: (H, D_{A'}) \ra (F(H), D_{B'})$
   implies that   for
$\epsilon=\la_1/\mu_1$:\newline
     (a)  $A'$ and
$\epsilon B'$ have the same real part Jordan form;\newline
  (b)  The   map
  $F|_H: (H, D^\epsilon_{A'}) \ra (F(H), D_{B'})$  is a $K$-quasisimilarity, where
  $K$ depends only on  $A'$, $B'$  and  $\eta$. \newline
     The statement (a) implies in particular $k_A-1=k_B-1$ (hence
     $k_A=k_B$),  $\la_i=\epsilon \mu_i$  and $e_i=d_i$  for $i<k_A$.
       Now it follows from Lemma   \ref{l4.6} that
       $\la_{k_A}=\epsilon\mu_{k_A}$.     If $A_{k_A}$ is a multiple of $I$, then
       (1)  follows from Lemmas \ref{l4.7}  and  \ref{l4.2}. If  $A_{k_A}$ is not a multiple of $I$,
  then Lemma \ref{l4.7} and Lemma \ref{l4.5} (1) imply that $A_{k_A}$  and $\epsilon B_{k_B}$
   have the same real part Jordan form.   Hence  (1) holds as well in this case.

 If $A_{k_A}$  is a multiple of $I$, then (2) follows from the statement
 (b) above and the arguments in the proof of Lemma \ref{l4.4}(2).
   If  $A_{k_A}$  is not a multiple of $I$, then Lemma \ref{l4.5}  (2)
  implies that $G: (V_{k}, D^\epsilon_{A_{k}})\ra (W_{k'}, D_{B_{k'}})$ is a  $K_1$-quasisimilarity
  with
   $K_1$ depending only on  $A_{k_A}$, $B_{k_B}$ and $\eta$.  %  and $\epsilon$.
 In this case, (2)  follows from this, (b) and the arguments in
   the proof of Lemma \ref{l4.4}(2).
% Section 4 of \cite{SX}.   When  applying the arguments in  Section 4 of
 %\cite{SX},
  % it is more convenient to use the quasimetric $\tilde D$, which is
   %biLipschitz equivalent to $D_A$.

\end{proof}

Next we will finish the proofs of the main theorems.
  % and their
   %corollaries.
  So  let $A, B$ be two  arbitrary  $n\times n$ matrices whose eigenvalues  have
positive  real parts. Let $G_A$, $G_B$ be equipped with arbitrary
admissible metrics.
 Then there are nonsingular matrices $P, Q$ such that
  $G_A$ is isometric to $G_{PAP^{-1}}$ (equipped with the standard metric) and
  $G_B$ is isometric to $G_{QBQ^{-1}}$  (equipped with the standard metric).
    Hence below in the proofs we will replace $(\R^n, D_A)$ and $(\R^n, D_B)$ with
    $(\R^n, D_{PAP^{-1}})$ and $(\R^n, D_{QBQ^{-1}})$  respectively.
  There also exist nonsingular matrices $P_0$ , $Q_0$ such that
  $G_{P_0AP_0^{-1}}$ and   $G_{Q_0BQ_0^{-1}}$ have pinched negative
  sectional curvature.
   We may choose  the same  $P_0AP_0^{-1}$ for all conjugate matrices $A$.
   Denote by $J$ and $J'$ the real part Jordan
  forms of $A$ and $B$ respectively.
By Proposition \ref{farb-m}, there are  biLipschitz maps
  $f_{J}:G_{P_0AP_0^{-1}} \ra G_J$  and
 $f_{P}:G_{P_0AP_0^{-1}} \ra  G_{PAP^{-1}}$.
   Then Corollary \ref{hou1}
       implies
  their  boundary maps $\p f_{J}: (\R^n, D_{P_0AP_0^{-1}})\ra (\R^n, D_J)$
  and  $\p f_{P}: (\R^n, D_{P_0AP_0^{-1}})\ra  (\R^n,
  D_{PAP^{-1}})$
    are also
  biLipschitz. Similarly, there are  biLipschitz maps
$f_{J'}:G_{Q_0BQ_0^{-1}} \ra G_{J'}$  and
 $f_{Q}:G_{Q_0BQ_0^{-1}} \ra  G_{QBQ^{-1}}$, whose
  boundary maps $\p f_{J'}: (\R^n, D_{Q_0BQ_0^{-1}})\ra (\R^n, D_{J'})$
  and  $\p f_{Q}: (\R^n, D_{Q_0BQ_0^{-1}})\ra  (\R^n,
  D_{QBQ^{-1}})$
    are also
  biLipschitz.

 %$f_Q:G_{QBQ^{-1}} \ra G_{J'}$  whose
  % boundary map $\p f_Q: (\R^n, D_{QBQ^{-1}})\ra (\R^n, D_{J'})$ is
  %also biLipschitz.

\noindent
  {\bf{Completing the proof of Theorem \ref{main1}}}.
 %Suppose $G_A$, $G_B$  are  equipped with arbitrary
%admissible metrics. Using the remarks above, we replace $G_A$ and
%$G_B$ with $G_{PAP^{-1}}$  and $G_{QBQ^{-1}}$ respectively.
 % First suppose $A$ and $sB$  have the same real part Jordan form.
  % Then $P_0AP_0^{-1}$ and $s Q_0BQ_0^{-1}$  also have the same real part Jordan form.
   %  By Proposition \ref{farb-m}, there is a biLipschitz map
    %  $f: G_{P_0AP_0^{-1}}\ra G_{Q_0BQ_0^{-1}}$.
  %Since $G_{P_0AP_0^{-1}}$  and $G_{Q_0BQ_0^{-1}}$ have pinched
 % negative sectional curvature, the boundary map
  % $\p f: (\R^n, D_{P_0AP_0^{-1}})\ra  (\R^n, D_{Q_0BQ_0^{-1}})$ is
   %quasisymmetric.
  %Since the two maps $\p f_P$ and $\p f_Q$ are biLipschitz, we see
  %that  $(\R^n, D_{PAP^{-1}})$  and $  (\R^n, D_{QBQ^{-1}})$
  %are quasisymmetric.
  The ``if" part follows from Proposition \ref{farb-m} since the
  boundary   map  of a quasiisometry between Gromov hyperbolic spaces is
  quasisymmetric.
 We will prove the  ``only if" part.
  So  we suppose $(\R^n, D_{PAP^{-1}})$  and $  (\R^n, D_{QBQ^{-1}})$
  are quasisymmetric.
Since the four  maps $\p f_P$, $\p f_J$,  $\p f_Q$ and $\p f_{J'}$
are biLipschitz, we see
  that  $(\R^n, D_{J})$  and $  (\R^n, D_{J'})$
  are quasisymmetric.
Now it    follows from  Lemma \ref{l4.2},
  Lemma \ref{l4.5}(1)  and Lemma \ref{l4.9}(1)
   that $J$ and $\epsilon J'$ have the same real part Jordan form,  where   $\epsilon=\la_1/\mu_1$.
     Hence
  $A$ and $\epsilon B$ also have the same real part Jordan form.

\qed

\b{Th}\label{main2text}
   Let $A$   and $B$ be   $n\times n$ matrices
whose eigenvalues all have positive real parts, and let $G_A$  and
$G_B$ be equipped  with arbitrary admissible metrics.
  Denote by  $\la_1$ and $\mu_1$  the smallest real parts of the eigenvalues
  of $A$ and $B$ respectively,  and   set $\epsilon=\la_1/\mu_1$.  %Let $k_A$ be the number of distinct
  %real parts of the eigenvalues of $A$.
  If  the real part Jordan form of $A$ is not a multiple of the identity
matrix $I_n$,  then  for every
  $\eta$-quasisymmetric map $F:(\R^n,  D_A)  \ra (\R^n, D_B)$,  the map
        $F:(\R^n,  D^\epsilon_A)  \ra (\R^n, D_B)$
    is  a  $K$-quasisimilarity, where $K$ depends only on  $\eta$,
     $A$, $B$ and the metrics on $G_A$, $G_B$.
  %       $\epsilon$ and $k_A$.

     \end{Th}

\noindent
  {\bf{Completing the proof of Theorem \ref{main2text}}}.
 Let  $F:  (\R^n, D_{PAP^{-1}})\ra   (\R^n, D_{QBQ^{-1}})$
  be an  $\eta$-quasisymmetric  map.
  Notice that the biLipschitz constant of the map $\p f_{J}$
  depends only on $A$ (actually the conjugacy class of $A$)  as the same $P_0AP_0^{-1}$ is chosen for all
  matrices  $A$ in the same conjugate class.
  However, the biLipschitz constant of $\p f_P$ depends on $P$ and
  hence on the admissible metric on $G_A$. Hence
   $\p f_J\circ \p f_P^{-1}$ is $L_1$-biLipschitz for some constant
   $L_1$  depending only on $A$ and the admissible metric on $G_A$.
   Similarly,  $\p f_{J'}\circ \p f_Q^{-1}$ is $L_2$-biLipschitz for some constant
   $L_2$  depending only on $B$ and the admissible metric on $G_B$.
It follows that
$$G:=(\p f_{J'}\circ \p f_Q^{-1})\circ F\circ (\p
f_J\circ \p f_P^{-1})^{-1}:  (\R^n, D_J)\ra (\R^n, D_{J'})$$
  is $\eta_1$-quasisymmetric,  where $\eta_1$ depends only on
   $L_1$, $L_2$  and $\eta$.
  Now  Lemma
\ref{l4.5}(2) and  Lemma \ref{l4.9}(2)  imply that
  $G$ is a  $K$-quasisimilarity, where $K$ depends only on  $J$, $J'$ and $\eta_1$.
 Consequently, $F$ is a $KL_1L_2$-quasisimilarity.

\qed

\section{Proof of the corollaries}\label{proof-co}

In this section we prove the corollaries from the introduction and
also derive a local version of Theorem \ref{main1}.

 Let $M$ be a  Hadamard manifold  with pinched negative
  sectional curvature,
$\xi_0\in \p M$, and  $x_0\in M$  a base point.  Let $\ga$ be the
geodesic with $\ga(0)=x_0$ and $\ga(\i)=\xi_0$.  Let $h_M=-B_\ga:
M\ra \R$, where $B_\ga$ is  the  Busemann function associated with
$\ga$. Set $H_t=h_M^{-1}(t)$.
  % The distance on the horosphere $H_t$ is denoted by $d_{H_t}$.
% Denote by $H_t$ the horosphere
 %  centered at $\xi_0$ and through $\ga(t)$.
   %\backslash \{\xi_0\}Let $\ga$ be the geodesic with $\ga(0)=x_0$ and $\ga(\i)=\xi_0$.
 A  parabolic visual quasimetric  $D_{\xi_0}$ on $\p
M\backslash\{\xi_0\}$ is defined as follows.
   For $\xi, \eta\in \p
M\backslash \{\xi_0\}$,
  $D_{\xi_0}(\xi, \eta)=e^t$ if and only if
   $\xi\xi_0\cap H_t$ and $\eta\xi_0\cap H_t$ have distance $1$ in
   the horosphere $H_t$.

Let $N$ be another Hadamard manifold  with pinched negative
  sectional
curvature, and $f:M \ra N$ a quasiisometry.  For any $\xi\in \p M$
and $x\in M$, we  set $\xi'=\p f(\xi)$ and $x'=f(x)$.   Let $\ga'$
be the geodesic with $\ga'(0)=x'_0$ and $\ga'(\i)=\xi'_0$.
 Set  $h_N=-B_{\ga'}$, where $B_{\ga'}$ is the Busemann function
  associated with $\ga'$. Denote $H'_t=h_N^{-1}(t)$.
     Let $D_{\xi'_0}$ be the parabolic
visual quasimetric on $\p N\backslash \{\xi'_0\}$ with respect to
the base point $x'_0$.
  % The distance function on

%The horosphere centered at $\xi'_0$  and through $\ga'(t)$  is
%denoted by $H'(t)$.

\b{Le}\label{3equiv} %Let $M, N$ be two Hadamard manifolds with
  %pinched negative sectional curvature, and $f:M\ra N$ a
  %quasiisometry.
 %Let $\xi_0\in \p M$ and set $\eta_0=\p f(\xi_0)$.
  Let $s>0$.
  Then the following three conditions are equivalent:\newline
  (1) there is a constant $C\ge 0$ such that the Hausdorff distance
    $HD(f(H_t), H'_{st})\le C$ for all $t$;\newline
  (2)   the boundary map $\p f: (\p M\backslash \{\xi_0\}, D^s_{\xi_0})\ra (\p N\backslash \{\xi'_0\},
  D_{\xi'_0})$  is biLipschitz;\newline  % where $D_{\xi_0}$ and $D_{\eta_0}$
     %are parabolic visual quasimetrics defined in ??;\newline
   (3)  there exists a  constant $C\ge 0$ such that
  $s\cdot d(x, y)-C\le d(f(x), f(y))\le s\cdot d(x,y)+C$     for all  $x,y\in M.$

   %$f$ is  an almost similarity.

\end{Le}

\b{proof} The arguments in the proof of  Lemma 6.4 in \cite{SX}
shows $(2)\Longrightarrow (1)$,    while the  arguments at the end
of \cite{SX} (proof of Corollary 1.2) yields $(1)\Longrightarrow
(3)$.
  We shall  prove $(3)\Longrightarrow (1)$  and
$(1)\Longrightarrow (2)$.

$(1)\Longrightarrow (2)$: Suppose (1) holds.
 Let $\xi\not=\eta\in \p M\backslash\{\xi_0\}$. Assume
  $D_{\xi_0}(\xi, \eta)=e^t$ and
$D_{\xi'_0}(\xi', \eta')=e^{t'}$. Let $\ga_\xi$ be  the geodesic
joining $\xi$ and $\xi_0$ with
    $\ga_\xi(0)\in H_0$ and $\ga_\xi(\i)=\xi_0$.
 By Lemma  6.2 of \cite{SX},  $\ga_\xi(t)$ is a $C_1$-quasicenter of
 $\xi, \eta, \xi_0$, and $\ga_{\xi'}(t')$ is  a $C_1$-quasicenter of
 $\xi'$, $\eta'$, $\xi'_0$, where $C_1$ depends only on
 the curvature bounds of $M$ and $N$.
        Since $f$ is a quasiisometry,
  $f(\ga_\xi(t))$ is a $C_2$-quasicenter
  of   $\xi'$, $\eta'$, $\xi'_0$, where $C_2$ depends only
  on $C_1$, the quasiisometry constants of $f$ and the curvature
  bounds of $N$.  It follows that $d(f(\ga_\xi(t)),\ga_{\xi'}(t'))\le C_3$, where $C_3$
   depends only on  $C_1$, $C_2$ and the curvature  bounds of $N$. By
   condition  (1), the point  $ f(\ga_\xi(t))$ is within $C$ of
  $H'_{st}$.  It follows that $\ga_{\xi'}(t')\in H'_{t'}$ is within $C+C_3$ of $H'_{st}$ and so
   $|t'-st|\le C+C_3$. Therefore,
   $e^{-(C+C_3)} e^{st}\le D_{\xi'_0}(\xi', \eta')=e^{t'}\le e^{C+C_3} e^{st}$.

$(3)\Longrightarrow (1)$:
 Suppose (3) holds.
 Let  $\omega: \R\ra M$ be any geodesic with $\omega(0)\in H_0$ and $\omega(\i)=\xi_0$.
   Then $f\circ \omega$ is a $(L_1,
C_1)$-quasigeodesic in $N$, where $L_1$ and $C_1$ depend
   only on $s$
and  $C$.  By the stability of quasigeodesics in a Gromov hyperbolic
space,  there is a constant $C_2$ depending only on $L_1, C_1$ and
the Gromov hyperbolicity constant of $N$, and a complete geodesic
$\omega'$ in $N$ with one endpoint $\xi'_0$  such that the Hausdorff
distance between $\omega'(\R)$ and $f\circ \omega (\R)$ is at most
$C_2$.
   %Hence $d(f(\omega(t)), \omega'(t))\le 2C_2$.
 Let  $t_1<t_2$.
  Then it follows from  condition (3) and the triangle inequality
  that
     $$|h_N(f(\omega(t_2)))-h_N(f(\omega(t_1)))-s(t_2-t_1)|\le   C_3,$$
   where $C_3$ depends only on $C$, $C_2$ and the Gromov
   hyperbolicity constant of $N$.
In particular, this applied to $\omega=\ga$, $t_2=t$ and $t_1=0$
  (or $t_2=0$ and $t_1=t$ if $t<0$)
implies  $|h_N(f(\ga(t)))-st|\le
  C_3$.
 %that $f(\ga(t))$ is within $C_3$ of $H'_{st}$.

Let $x\in H_t$ be arbitrary.  Let $\omega_1$ be the geodesic with
$\omega_1(t)=x$ and $\omega_1(\i)=\xi_0$.  Pick any $t_2\ge t$ with
 $d(\ga(t_2), \omega_1(t_2))\le 1$.
  By condition (3),
  $$|h_N(f(\ga(t_2)))-h_N(f(\omega_1(t_2)))|\le   d(f(\ga(t_2)), f(\omega_1(t_2)))\le s +C.$$
  The discussion from the
 preceding paragraph implies
$$|h_N(f(\omega_1(t_2)))-h_N(f(\omega_1(t)))-s(t_2-t)|\le
  C_3 $$
and
   $$|h_N(f(\ga(t_2)))-h_N(f(\ga(t)))-s(t_2-t)|\le
  C_3.$$
 These inequalities together with the one at the end of last
 paragraph imply
   $$|h_N(f(\omega_1(t)))-st|\le C_4:=3C_3+s+C.$$
  Hence $f(x)=f(\omega_1(t))$ is within $C_4$ of $H'_{st}$.
   This shows
 $f(H_{t})$ lies in the
$C_4$-neighborhood of  $H'_{st}$. By considering a quasi-inverse of
$f$, we see that the Hausdorff
 distance
  $HD(f(H_t), H'_{st})\le C_5$, where $C_5$ depends only on $s$, $C$
  and the Gromov hyperbolicity constants of $M$ and $N$.

% between them is at most $C'$ for some constant $C'$
% depending only on $C$ and the Gromov hyperbolicity constants of $M$
 %and $N$.  The inequality at the end of the preceding paragraph
 %applied to $\ga_1$, $t_1=t, t_2=0$ implies
  %$|h_N(\ga_1(t))-t|\le C+4C_2+|h_N(\ga_1(0))|$.
%Therefore, the Hausdorff distance from $f(H_t)$ to ?? is at most
 %$C'+C+4C_2+|h_N(\ga_1(0))|$.

\end{proof}

A local version of Theorem \ref{main1}  also holds:

\b{Th}\label{localtheorem1}  Let $A$, $B$ be $n\times n$ matrices
 whose eigenvalues have positive real parts,  and let $G_A$ and
 $G_B$ be equipped with arbitrary admissible metrics.
   Let $U\subset (\R^n, D_A)$, $V\subset (\R^n,
D_B)$ be open subsets, and
  $F: (U, D_A)\ra (V, D_B)$ an  $\eta$-quasisymmetric map.
     Then $A$   and
$s B$  have the same real part Jordan form
  for some $s>0$.
%     where
  %$\epsilon=\la_1/\mu_1$.

\end{Th}

\b{proof}
  By Corollary \ref{hou1} and the discussion before the proof of Theorem \ref{main1}
        we may  assume $A$ and $B$ are in real
  part Jordan form.  Fix a base point $x\in U$. We may assume  both $x$ and
$F(x)$ are the origin $o$.
 Then there is some constant $a>1$ and a
sequence of  distinct triples $(x_k, y_k, z_k)$ from $U$ satisfying
$x_k=o$, $D_A(x_k, y_k)\ra 0$ and
$$\frac{D_A(x_k, y_k)}{D_A(x_k, z_k)},\frac{D_A(y_k, x_k)}{D_A(y_k, z_k)},\frac{D_A(z_k, x_k)}{D_A(z_k,
y_k)}\in ({1}/{a}, a).$$
 Such a triple   can be chosen from the  eigenspace of $\la_1$
 (the smallest eigenvalue of $A$)  so
 that $x_k=o$ is the middle point of the segment $y_kz_k$.
  Since $F$ is $\eta$-quasisymmetric,  there is
a  constant $b>0$ depending only on $a$ and $\eta$ such that:
$$\frac{D_B(F(x_k), F(y_k))}{D_B(F(x_k), F(z_k))},\frac{D_B(F(y_k), F(x_k))}{D_B(F(y_k), F(z_k))},
  \frac{D_B(F(z_k), F(x_k))}{D_B(F(z_k),  F(y_k))}\in ({1}/{b}, b).$$
  Assume $D_A(x_k, y_k)=e^{-t_k}$ and $D_B(F(x_k),
  F(y_k))=e^{-t'_k}$.  Then
$e^{t_kA}: (U,  e^{t_k}D_A)\ra (e^{t_k A}U, D_A)$ is an isometry.
Hence the sequence of pointed metric spaces $(U,e^{t_k}D_A, o)$
converges (as $k \ra \i$) in the pointed Gromov-Hausdorff topology
  towards $(\R^n, D_A)$.  Similarly,
 the sequence of pointed metric spaces $(V,  e^{t'_k}D_B, o)$
converges (as $k \ra \i$) in the pointed Gromov-Hausdorff topology
  towards $(\R^n, D_B)$.
  On the other hand,   the sequence of maps $F_k=F:(U,  e^{t_k}D_A)\ra (V,
        e^{t'_k}D_B)$ are all $\eta$-quasisymmetric, and
          %the above separation conditions still hold for
         the triples
         $(x_k, y_k, z_k)\in  (U,  e^{t_k}D_A)$  and  $(F(x_k), F(y_k),
         F(z_k))\in (V,
        e^{t'_k}D_B)$  are uniformly separated and uniformly bounded.
    Now the  compactness  property of quasisymmetric maps implies   that a subsequence  of $\{F_k\}$ converges in the
  pointed Gromov-Hausdorff  topology  towards an
  $\eta$-quasisymmetric map  $F': (\R^n, D_A)\ra (\R^n, D_B)$.
  Now the theorem follows from Theorem \ref{main1}.

  % On the other hand, since $e^{t_kA}: (U,  e^{t_k}D_A)

\end{proof}

\b{Le}\label{ati} Let $F: \p G_A\ra \p G_B$ be a quasisymmetric map,
where $\p G_A$ and $\p G_B$ are equipped with visual metrics.
   Let $\xi_0\in \p G_A, \xi'_0\in \p G_B$  be the points
   corresponding to upward oriented vertical
  geodesic rays.  If the real part Jordan form of $A$ is not a
  multiple of  the identity matrix, then
   $F(\xi_0)=\xi'_0$.

\end{Le}

\b{proof}  The proof is similar to that of Proposition 3.5,
\cite{X}.
 Suppose  $F(\xi_0)\not=\xi'_0$.
  By the relation
 between visual metrics and parabolic visual metrics (see
 Section 5 of \cite{SX}), the map
  $$F: (\R^n \backslash \{F^{-1}(\xi'_0)\}, D_A)\ra (\R^n\backslash
  \{F(\xi_0),  D_B)$$
     is locally quasisymmetric.  By  Theorem  \ref{localtheorem1},
   $A$ and $s B$ have the same real part Jordan form for some
   $s>0$.  In particular, we have $k_B=k_A$;
    the fibers of $\pi_A$ and $\pi_B$ have the same dimension if $k_A=1$,
       and the subspaces $\prod_{i<k_A}V_i$ and $\prod_{j<k_B}W_j$ have the same dimension if $k_A\ge 2$.
   If $k_A=1$, let $H$ be a fiber of $\pi_A$  not containing
   $F^{-1}(\xi'_0)$; if $k_A\ge 2$, then let $H$ be an affine
    subspace parallel to $\prod_{i<k_A}V_i$ and not containing
$F^{-1}(\xi'_0)$.
 Let $m$ be the topological dimension of $H$.
 Then $H\cup \{\xi_0\}\subset \p G_A$ is an $m$-dimensional topological sphere.
  Since $F(\xi_0)\not=\xi'_0$ and   $F^{-1}(\xi'_0)\notin H$,
   the image $F(H\cup \{\xi_0\})$ is a $m$-dimensional topological
   sphere in $\R^n=\p G_B\backslash \{\xi'_0\}$.
 In particular, $F(H\cup \{\xi_0\})$ (and hence $F(H)$)  is not contained in any fiber
 of $\pi_B$ (if $k_A=1$) or any affine subspace parallel to
  $\prod_{j<k_B}W_j$ (if $k_A\ge 2$).
  Now the arguments of Lemma \ref{l4.3}   and Lemma \ref{l4.7}  yield a contradiction.
  Hence   $F(\xi_0)=\xi'_0$.

\end{proof}

%\section{Proofs of the Corollaries}\label{proof-co}

Now Corollary \ref{c1}   follows from Proposition \ref{farb-m},
 Lemma \ref{ati},  Theorem \ref{main1}  and the
 fact that a   quasiisometry between Gromov hyperbolic spaces
 induces a quasisymmetric map  between the ideal boundaries.

\noindent
  {\bf{Proofs of  Corollary \ref{c2} and Corollary  \ref{c3}}}.
    We use the notation introduced before the proof of Theorem
    \ref{main1}.
  Let $f:  G_{PAP^{-1}}  \ra   G_{QBQ^{-1}}$
  be a quasiisometry.  By Lemma \ref{ati}, $f$ induces  a
    boundary map $\p f: (\R^n,  D_{PAP^{-1}})\ra (\R^n,
   D_{QBQ^{-1}})$, which  is quasisymmetric.
  By Theorem \ref{main2},  there is some $s>0$ such that
    $\p f: (\R^n,  D^s_{PAP^{-1}})\ra (\R^n,
   D_{QBQ^{-1}})$
    is biLipschitz.  Since $\p f_P$ and $\p f_Q$ are also
   biLipschitz,  the boundary map
        $\p (f^{-1}_Q\circ f \circ f_P):    (\R^n,  D^s_{P_0AP_0^{-1}})\ra (\R^n,
   D_{Q_0BQ_0^{-1}}) $
      of $f^{-1}_Q\circ f \circ f_P:
      G_{P_0AP_0^{-1}}\ra  G_{Q_0BQ_0^{-1}}$  is biLipschitz. Since
  $G_{P_0AP_0^{-1}}$  and  $G_{Q_0BQ_0^{-1}}$
  have pinched negative sectional curvature, Lemma   \ref{3equiv}
        implies
   the map  $f^{-1}_Q\circ f \circ f_P$ is height-respecting and is
   an almost similarity. By Proposition \ref{farb-m}  and
       Corollary \ref{hou1}
        the two maps $f_P$ and $f_Q$ are
     height-respecting and are
    almost  similarities.
    Hence $f$ is   height-respecting and is
   an almost similarity.
  \qed

The proof of Corollary \ref{c4} is the same as in \cite{SX}
(Corollary 1.3).

Next we give  a proof of Corollary \ref{c5}.   Recall that a  group
$G$ of quasisimilarity  maps of $(\R^n, D_A)$ is a uniform group if
there is some $K\ge 1$ such that every element of $G$ is a
$K$-quasisimilarity. Dymarz  and Peng have established the following
(see \cite{DP} for  the definition of almost homotheties):

\b{theorem}\label{d-p} \e{(\cite{DP})}
  Let $A$ be a  square matrix whose eigenvalues
all have positive real parts, and $G$ be a
   uniform
group of  quasisimilarity  maps of $(\R^n, D_A)$.
  If the induced action of $G$ on the space of distinct couples of
    $\R^n$  is cocompact, then $G$
 can be conjugated
by a biLipschitz map into the
 group of almost homotheties.

   \end{theorem}

\noindent
 {\bf{Proof of Corollary \ref{c5}.}}
  Let $G$ be a group of quasim\"obius maps of $(\p G_A, d)$ such that
  every element of $G$ is $\eta$-quasim\"obius, where $d$ is a fixed
  visual metric on $\p G_A$.
 Let  $\xi_0\in \p G_A$
be the point corresponding to vertical geodesic rays.  Since the
real part Jordan form of $A$
  is not a multiple of the identity matrix,
% certain foliations of
 % $\R^n$ by parallel affine subspaces are preserved by all quasisymmetric maps $\p G_A\ra \p G_A$
  %(see the proofs of Theorem \ref{main1}  and Theorem \ref{main2},
  %in particular, Lemma \ref{l4.3} and Lemma \ref{l4.7}).
   %  Since the closures (in $\p G_A=\R^n\cup \{\xi_0\}$)  of two distinct leaves intersect in a single
    % point $\{\xi_0\}$,
      Lemma \ref{ati} implies that
     the point $\xi_0$ is fixed by all
     quasisymmetric maps $\p G_A\ra \p G_A$.
       Hence  $G$ restricts to a group of quasisymmetric maps of
       $(\R ^n, D_A)$.
For any three distinct points $\xi_1, \xi_2, \xi_3\in \R^n=\p
G_A\backslash \{\xi_0\}$, the quasim\"obius condition applied to
 the quadruple $Q=(\xi_1, \xi_2, \xi_3, \xi_0)$ implies that
every element of $G$ is an $\eta$-quasisymmetric map of
 $(\R ^n, D_A)$.  Now Theorem \ref{main2} implies that there is some
 $K\ge 1$ such that every element of $G$ is  a $K$-quasisimilarity.
 In other words, $G$ is a uniform group of quasisimilarities
 of   $(\R ^n, D_A)$.

 Since the induced action of $G$ on the space of distinct triples of
 $(\p G_A, d)$ is cocompact, there is some $\delta>0$ such that for
 any distinct triple $(\xi_1, \xi_2, \xi_3)$, there is some $g\in G$
 such that  $d(g(\xi_i), g(\xi_j))\ge \delta$ for all $1\le i\not= j\le
 3$. Now let $\xi\not= \xi_2\in \R^n=\p G_A\backslash\{\xi_0\}$ be
 any distinct couple.  Then there is an element $g\in G$ as above corresponding to  the triple $(\xi_0, \xi_1,
 \xi_2)$.   Since  $g(\xi_0)=\xi_0$,  there are two constants $a,
 b>0$ depending only on $\delta$ such that $D_A(g(\xi_1), o)\le b$,
   $D_A(g(\xi_2), o)\le b$ and  $D_A(g(\xi_1), g(\xi_2))\ge a$.  This
 shows that $G$ acts cocompactly on the space of distinct couples of
 $(\R^n, D_A)$.

   Now the corollary follows from the theorem of Dymarz-Peng.

\qed

\section{QS maps in the Jordan block case}\label{jor}

In this section we describe all the quasisymmetric maps on the ideal
boundary in the case when $A$ is a Jordan block.
  % As a consequence,
%we are able to establish a Liouville type theorem for $G_A$, see
  %Theorem
   %\ref{liu}.

%In this section we characterize all quasisymmetric maps in the
%Jordan block case, prove a Liouville type theorem and improve
%Corollary ??

\b{Th}\label{Jordan} Let $J_n=I_n +N$ be  the   $n\times n$ ($n\ge
2$) Jordan
 block with  eigenvalue $1$.  Then a bijection
 $F:(\R^n, D_{J_n})\ra (\R^n, D_{J_n})$
   is a quasisymmetric map   if  and only of
 there are constants $a_0\not=0, a_1, \cdots, a_{n-2}\in \R$, a vector
       $v\in \R^n$  and a
 Lipschitz map $C:\R\ra \R$ such that
   $$F(x)=(a_0 I_n +a_1 N +\cdots +a_{n-2} N^{n-2})x+ v+ \tilde{C}(x)$$
       for
   all $x=(x_1, \cdots, x_n)^T\in \R^n$, where  $\tilde{C}(x)=(C(x_n), 0, \cdots, 0)^T$. Here $T$
   indicates matrix transpose.

\end{Th}

We first prove that  every map of the indicated form is  actually
  biLipschitz.  Notice that the map $F$ described in the theorem
  decomposes  as
    $F=F_1\circ F_2\circ F_3$, with
    $F_1(x)=x+v$,   $F_2(x)=x+\tilde{C_1}(x)$
    and $F_3(x)=(a_0 I_n +a_1 N +\cdots +a_{n-2}
     N^{n-2})x$,  where $C_1:\R\ra \R$  is defined by  $C_1(t)=C(t/a_0)$.   Since  $D_{J_n}$ is invariant under Euclidean
     translations,  $F_1$ is an isometry.  We shall prove that $F_2$
     and $F_3$ are biLipschitz in the next two lemmas.

   For  an
 $n\times n$ matrix $M=(m_{ij})$,
     set  $Q(M)=\sum_{i,j}
 m_{ij}^2$.
   % The operator norm  of $M$ is
 %$||M||=\sup\{|Mx|:   x\in \R^n, |x|=1\}$.
   % Let $Q(M)=\sum_{i,j}
   %m_{ij}^2$.
 We will use the fact
  $||M||\le  Q(M)^{\frac{1}{2}}$, where $||M||$ denotes the operator norm of $M$.     %(\sum_{i,j} m_{ij}^2)^{\frac{1}{2}}$.

\b{Le}\label{f1} Suppose  $C:\R\ra \R$ is $L$-Lipschitz  for some
$L>0$. Then $F_2: (\R^n, D_{J_n})\ra (\R^n, D_{J_n})$, $F_2(x)=x+
\tilde C(x)$ is $L'$-biLipschitz, where $L'$ depends only on $L$ and
the dimension $n$.

\end{Le}

\b{proof}%Assume $C:\R\ra \R$ is  $L$-Lipschitz  for some $L\ge
    %1$.
 Let $x=(x_1, \cdots, x_n)^T$ and $x'=(x'_1, \cdots, x'_n)^T$ be two
 arbitrary points in $\R^n$.   Then $F_2(x)=(x_1+ C(x_n), x_2,
 \cdots, x_n)^T$  and $F_2(x')=(x'_1+C(x'_n), x'_2, \cdots,
 x'_n)^T$.
   Assume  $D_{J_n}(x,x')=e^t$ and
 $D_{J_n}(F_2(x), F_2(x'))=e^s$.
 We need to show that there is some constant  $a$ depending only on
 $L$ and $n$ such that $|t-s|\le a$.

Since $D_{J_n}(x,x')=e^t$,   we  have  %the definition of $D_A$ implies
  %$|e^{-tA}(x'-x)|=1$, and hence
 $e^t=|e^{-tN}(x'-x)|$, see Section \ref{metriconb}. Similarly,
   $D_{J_n}(F_2(x), F_2(x'))=e^s$ implies
    $e^s=|e^{-sN}(F_2(x')-F_2(x))|$.
   Notice $F_2(x')-F_2(x)=(x'-x)+ w$, where $w= (C(x'_n)-C(x_n), 0, \cdots, 0)^T$.
      The only nonzero entry in $e^{-tN}w$ is
   $C(x'_n)-C(x_n)$. So we have
     $$|e^{-tN}w|=
   |C(x'_n)-C(x_n)|\le L |x'_n-x_n|.$$
   On the other hand,
 the last entry in
    $e^{-tN}(x'-x)$  is $(x'_n-x_n)$, hence
       $$|e^{-tN}w|\le  L |x'_n-x_n|\le L |e^{-tN}(x'-x)|=Le^t.$$
We write
  $$e^{-sN}(F_2(x')-F_2(x))=e^{(t-s)N}e^{-tN}[(x'-x)+w]=e^{(t-s)N}[e^{-tN}(x'-x)+
  e^{-tN}w].$$
                Now
\begin{align*}
e^s   & =|e^{-sN}(F_2(x')-F_2(x))|=|e^{(t-s)N}[e^{-tN}(x'-x)+
  e^{-tN}w]|  \\
   & \le ||e^{(t-s)N}|| \cdot |e^{-tN}(x'-x)+
  e^{-tN}w|\le ||e^{(t-s)N}|| \cdot \left\{|e^{-tN}(x'-x)|+
   |e^{-tN}w|\right\}\\
   & \le||e^{(t-s)N}||\cdot \big\{e^t+L e^t\big\}\le
   e^t(1+L)\sqrt{Q(e^{(t-s)N})}.
\end{align*}
From this we derive
  $e^{s-t}\le (1+L)\sqrt{Q(e^{(t-s)N})}$.
   Notice that $Q(e^{(t-s)N})$  is a
     polynomial of degree $2(n-1)$ in $t-s$ that depends only on
 $n$.
   It follows that there is a  constant $a$ depending only on $n$ and
   $L$ such that $s-t\le a$.
     Since the inverse of $F_2$ is $F_2^{-1}(x)=x+(-C(x_n), 0,
     \cdots, 0)^T$, the above argument applied to $F_2^{-1}$ yields
     $t-s\le a$. Hence $|s-t|\le a$, and we are done.

\end{proof}

\b{Le}\label{f3}  Let  $F_3:(\R^n, D_{J_n})\ra (\R^n, D_{J_n})$ be
  given
by
$$F_3(x)= (a_0 I_n +a_1 N +\cdots +a_{n-1}
     N^{n-1})x,$$
         where
   $a_0\not=0, a_1, \cdots, a_{n-1}\in\R$ are
constants. Then $F_3$   is  $L$-biLipschitz for some $L$ depending
only on $n$ and $a_0, a_1, \cdots, a_{n-1}$.

\end{Le}

\b{proof}
  The proof is similar to that of  Lemma  \ref{f1}.
 Let $x, x'\in \R^n$ be arbitrary. Assume $D_{J_n}(x, x')=e^t$ and
 $D_{J_n}(F_3(x), F_3(x'))=e^s$.
   Then we have
    $e^t=|e^{-tN}(x'-x)|$  and
     $e^s=|e^{-sN} (F_3(x')-F_3(x))|$.  We need to find a constant
     $a$ that depends only on $n$ and the numbers $a_0$, $\cdots$,
     $a_{n-1}$ such that $|s-t|\le a$.

       Set $B_1=e^{(t-s)N}$ and $B_2=a_0 I_n
+a_1 N +\cdots +a_{n-1}
     N^{n-1}$.  Notice that $B_2$ commutes with $N$.
       We
        have
       \begin{align*}
       e^s& =|e^{-sN} (F_3(x')-F_3(x))|=|e^{(t-s)N} e^{-tN}
       B_2(x'-x)|\\
       &=|B_1  B_2  e^{-tN}(x'-x)|
       \le ||B_1||\cdot ||B_2||\cdot |e^{-tN}(x'-x)|\\
       &\le  \sqrt{Q(B_1)} \sqrt{Q(B_2)}\;e^t.
\end{align*}
 Hence  $e^{s-t}\le \sqrt{Q(B_1)Q(B_2)}$.     Since
   $Q(B_1) Q(B_2)$   is a  polynomial in $t-s$ that  depends only on $n$ and the
   numbers $a_0, \cdots, a_{n-1}$,
          there is some constant  $a>0$ depending only on $n$ and $a_0, \cdots, a_{n-1}$
             such that
   $s-t\le a$.

Notice that $F_3^{-1}(x)=B_2^{-1}x$.
  Set
  $$\beta=-\left(\frac{a_1}{a_0}N+\cdots+ \frac{a_{n-1}}{a_0}N^{n-1}\right).$$
    Then $\beta^n=0$.
  We have $B_2=a_0(I-\beta)$
and $B_2^{-1}=a_0^{-1}(I+\beta+\beta^2+\cdots +\beta^{n-1})$.
  It follows that  $B_2^{-1}$ has the expression
 $B_2^{-1}=a_0^{-1} I+ b_1N+\cdots +b_{n-2}N^{n-2}+b_{n-1}N^{n-1}$, where $b_1,
 \cdots, b_{n-1}$ are constants depending only on $a_0, \cdots, a_{n-1}$.
   Now the preceding paragraph implies that
   $t-s\le a'$ for some constant  $a'$
       depending only on $n$ and $a_0^{-1}, b_1,
   \cdots, b_{n-1}$, hence only on  $n$ and  $a_0, \cdots, a_{n-1}$.
     Therefore  $|s-t|\le \max\{a, a'\}$, and the proof of Lemma \ref{f3}
    is complete.

\end{proof}

To prove that every quasisymmetric map has the described type, we
induct on  $n$. The basic step $n=2$ is given by Theorem \ref{xthm}.
 Now we assume $n \ge 3$ and that Theorem \ref{Jordan}
     holds for  $J_{n-1}$.

Let  $F: (\R^n, D_{J_n})\ra (\R^n, D_{J_n})$ be a
     quasisymmetric  map.
 Let $\mathcal {F}_i$ ($i=1, \cdots, n-1$) be
 the foliation of $\R^n$ consisting of affine subspaces parallel to
 the  linear subspace
  $$H_i:=\{x=(x_1, \cdots, x_n)^T\in \R^n:
 x_{i+1}=\cdots=x_n=0\}.$$
      Then the proof  of    Theorem \ref{main2}
  shows that
  the foliation $\mathcal {F}_i$  is preserved by  $F$.
    %quasisymmetric map of $(\R^n, D_{J_n})$.
  To be more precise, if
     $H$ is an affine subspace parallel to
     $H_i$,  %  and if $F: (\R^n, D_{J_n})\ra (\R^n, D_{J_n})$ is
         %$\eta$-quasisymmetric,
      then $F(H)$ is also an affine subspace parallel
     to $H_i$.   In particular,
  $F$ maps every line parallel to the  $x_1$-axis  (that is, parallel to $H_1$) to a line parallel to  the $x_1$-axis,
  and maps every horizontal hyperplane (that is, parallel to $H_{n-1}$)  to a horizontal hyperplane.
It follows that there is a map $G: \R^{n-1}
 \ra    \R^{n-1}$ such that for any $y\in \R^{n-1}$, $F(\R\times
\{y\})=\R\times \{G(y)\}$. For each $y\in \R^{n-1}$, there is a map
$H(\cdot, y): \R\ra \R$ such that $F(x_1,y)=(H(x_1, y), G(y))$.

Arguments similar to the proofs of Lemmas \ref{l11} and \ref{l12}
show the following:\newline
   (1)   for each $y\in \R^{n-1}$, the restriction of $D_{J_n}$
    to  $\R\times \{y\}$  agrees with  the Euclidean distance on
    $\R$;\newline
        (2) for any two $y_1, y_2\in \R^{n-1}$,  the Hausdorff
        distance with respect to $D_{J_n}$:
          $HD(\R\times \{y_1\},\R\times
          \{y_2\})=D_{J_{n-1}}(y_1, y_2)$;\newline
            (3) for any $p=(x_1, y_1)\in \R\times \R^{n-1}$ and any
            $y_2\in \R^{n-1}$,
               we have $D_{J_n}(p, \R\times \{y_2\})=D_{J_{n-1}}(y_1,
               y_2)$.\newline
 Hence  each
$H(\cdot, y): (\R, |\cdot|)  \ra (\R, |\cdot|)$
  is    quasisymmetric, and
  the arguments on page 10 of
 \cite{X}
   shows
that     $G: (\R^{n-1},  D_{J_{n-1}})
 \ra (\R^{n-1}, D_{J_{n-1}})$   is  also
       quasisymmetric.

Now the induction hypothesis applied to  $G$ shows that there are
constants $a_0\not=0$, $a_1$, $\cdots, a_{n-3}$, $b_i$ ($2\le i\le
n)$ and a Lipschitz map $g: \R\ra \R$ such that
 \[
 G\left(\begin{array}{c}x_2\\  x_3\\ \cdots\\ x_n\end{array}\right)=\left(\begin{array}{c}a_0 x_2+a_1x_3+\cdots +
 a_{n-3}x_{n-1}+b_2+g(x_n)\\ a_0x_3+a_1x_4+\cdots +a_{n-3}x_n+b_3\\
 \cdots\\
 a_0x_{n-1}+a_1x_n+b_{n-1}\\
 a_0x_n+b_n\end{array}\right).
 \]
  Notice that the horizontal hyperplane $\R^{n-1}\times \{x_n\}$  at height $x_n$ is mapped by
   $F$ to the horizontal hyperplane
$\R^{n-1}\times \{a_0x_n+b_n\}$
   at height $a_0x_n+b_n$.
      Since  the restriction of $D_{J_n}$ to a horizontal hyperplane
      agrees with  $D_{J_{n-1}}$ (Lemma \ref{l11}), the map
      $$F: (\R^{n-1}\times \{x_n\}, D_{J_{n-1}}) \ra (\R^{n-1}\times
      \{a_0x_n+b_n\},D_{J_{n-1}}) $$
          is    quasisymmetric. Now
      the induction hypothesis, the fact $F(x_1, y)=(H(x_1, y),
      G(y))$ and the expression of $G$ imply
      that
      $$H(x_1, y)=a_0x_1+a_1x_2+\cdots
      +a_{n-3}x_{n-2}+c_1(x_n)+c_2(x_{n-1}, x_n),$$
          where $c_1: \R\ra
      \R$ and $c_2:   \R^2\ra \R$ are two maps and for each fixed $v$,  $c_2(u, v)$ is
      Lipschitz in $u$.   Since $F$ is a homeomorphism,
   $c_1$  and $c_2$ are continuous.
      Define $c_3:\R^2\ra \R$ by $c_3(u,v)=c_1(v)+c_2(u, v)$.  After composing $F$ with a map
      of the described type,  we may assume $F$ has the following
      form
       $$F(x_1, x_2, \cdots, x_n)=(x_1+c_3(x_{n-1}, x_n),
       x_2+g(x_n), x_3, \cdots, x_n).$$
         We need to show that
  there are constants $a_{n-2}$,  $d_2$ and a Lipschitz map
  $C:\R\ra \R$ such that
   $g(x_n)=a_{n-2}x_n+d_2$  and $c_3(x_{n-1}, x_n)=a_{n-2}
   x_{n-1} +C(x_n)$.

%\end{document}

\b{Le}\label{l5.1}
 There is a constant $L$    %depending only on $L$ and ??
    such that the
 following holds  for all $u, v, v'\in \R$:
 $$\Big\vert\big\{c_3(u+(v'-v)\ln |v'-v|, v')-c_3(u,v)\big\}-\ln{|v'-v|}\big\{g(v')-g(v)\big\}\Big\vert \le L |v'-v|.$$

 \end{Le}

 \b{proof}
Let  $u, v, v'\in \R$.  Let $x\in \R^n$ with $x_{n-1}=u$,  $x_n=v$.
 Set    $t=\ln |v'-v|$  and
   let  $y=(y_1, \cdots, y_n)^T$ be the unique solution of
  $e^{-tN}y=(0, \cdots, 0,  v'-v)^T$.   Let $x'=x+y$.
Notice $y_n=v'-v$,  $\,y_{n-1}=(v'-v)\ln|v'-v|$,  $\,x'_n=v'$   and
  $$x'_{n-1}=x_{n-1}+y_{n-1}=u+(v'-v)\ln|v'-v|.$$
         Notice also that   $t$ is the
       smallest  solution for
       $e^t=|e^{-tN}(x'-x)|$ and so $D_{J_n}(x, x')=e^t$.
        Suppose $D_{J_n}(F(x), F(x'))=e^s$.  Then
          $e^s=\big\vert e^{-sN}(F(x')-F(x))\big\vert$.
           By Theorem \ref{main2}, $F$ is $L_1$-biLipschitz  for some $L_1\ge 1$. Hence
   $e^t/{L_1}\le e^s\le L_1 e^t$.  It follows that  $|t-s|\le \ln {L_1}$.
    Now we write
\b{align*} e^{-sN}(F(x')-F(x))& =e^{-sN}(x'-x)+
      e^{-sN}\left({\begin{array}{c}c_3(x'_{n-1},
      x'_n)-c_3(x_{n-1}, x_n)\\  g(x'_n)-g(x_n)\\ 0\\ \cdot \\ \cdot \\ \cdot\\
      0\end{array}}\right)\\
  &=e^{(t-s)N}e^{-tN}(x'-x)+
  e^{(t-s)N}e^{-tN}\left({\begin{array}{c}c_3(x'_{n-1},
      x'_n)-c_3(x_{n-1}, x_n)\\  g(x'_n)-g(x_n)\\ 0\\ \cdot \\ \cdot \\ \cdot \\
      0\end{array}}\right)\\
  &
      =e^{(t-s)N}\left(\begin{array}{c}\big\{c_3(x'_{n-1},
      x'_n)-c_3(x_{n-1}, x_n)\big\}-t\big\{g(x'_n)-g(x_n)\big\}\\  g(x'_n)-g(x_n)\\ 0\\ \cdot \\ \cdot \\ \cdot \\
      0\\
      x'_n-x_n\end{array}\right).
\end{align*}
   Set
    $$\tau=\big\{c_3(x'_{n-1},
      x'_n)-c_3(x_{n-1}, x_n)\big\}-t\big\{g(x'_n)-g(x_n)\big\}.$$
        The first entry of $e^{-sN}(F(x')-F(x))$ is
        $$q:=\tau
      +(t-s)\big\{g(x'_n)-g(x_n)\big\}+\frac{(t-s)^{n-1}}{(n-1)!}(x'_n-x_n).$$
         %Denote this term by $q$.
       We have
       $$|q|\le
       |e^{-sN}(F(x')-F(x))|=e^s\le L_1 e^t=L_1|v'-v|.$$
          Recall that $g$ is $L_2$-Lipschitz for some $L_2\ge 0$.
           Hence, $$|g(x'_n)-g(x_n)|\le
       L_2|x'_n-x_n|=L_2|v'-v|.$$
         Now it follows from
        $|t-s|\le \ln {L_1}$ and
           the triangle inequality that
$$|\tau|\le \left(L_1+L_2\ln
      {L_1}+\frac{(\ln {L_1})^{n-1}}{(n-1)!}\right)|v'-v|.$$

 \end{proof}

Recall that the map $g$
 is Lipschitz  and for each fixed $v$,  $c_3(u, v)$ is Lipschitz in $u$.
 Hence $g$ is differentiable a.e., and for each fixed $v$, the
 partial derivative $\frac{\partial c_3}{\partial u}$ exists
 for a.e. $u$.

 \b{Le}\label{l5.2}
Let $v$ be any point such that $g'(v)$ exists. Then
 $c_3(u,v)=c_3(0, v)+g'(v) u$ for all $u$.
 % $c_3(x_{n-1}, v)$ is affine in $x_{n-1}$.

 \end{Le}

 \b{proof}
Fix an arbitrary $u\in \R$.
 Let $a>0$.   For any positive integer $n$,   define $(y_0, z_0)=(u,v)$ and
  $(y_i, z_i)=(u+i\frac{a}{n}\ln \frac{a}{n}, v+i\frac{a}{n})$
  ($1\le i\le n$).   Applying Lemma \ref{l5.1}  to
   $y_{i-1}, z_{i-1},  z_i$
  we obtain:
  $$\Big\vert\big\{c_3(y_i, z_i)-c_3(y_{i-1}, z_{i-1})\big\}-\ln
  \frac{a}{n}\big\{g(z_i)-g(z_{i-1})\big\}\Big\vert \le L\frac{a}{n}.$$
  Now let $k=k(n)$ be the integer part of $n/{\ln\frac{n}{a}}$. Then
    $\frac{k}{n}\ln \frac{a}{n}\ra -1$ as $n\ra \i$.
      Combining the above inequalities for $1\le i\le k$ and using
         the
      triangle inequality, we obtain
$$\Big\vert\big\{c_3(y_k, z_k)-c_3(u, v)\big\}-\ln
  \frac{a}{n}\big\{g(z_k)-g(v)\big\}\Big\vert\le L\frac{ak}{n}.$$
Now divide   both sides by $\frac{ak}{n}\ln\frac{n}{a}$ (which
converges to $a$ as $n\ra \i$), we get:
$$\bigg\vert\frac{\big\{c_3(y_k, z_k)-c_3(u, v)\big\}}{\frac{ak}{n}\ln\frac{n}{a}}
+
  \frac{\big\{g(z_k)-g(v)\big\}}{\frac{ak}{n}}\bigg\vert\le \frac{L}{\ln \frac{n}{a}}.$$
  As $n\ra \i$, we have $z_k=v+\frac{ak}{n}\ra v$, $y_k\ra u-a$.  Also, since
  $g'(v)$ exists, we have
  $$\frac{\big\{g(z_k)-g(v)\big\}}{\frac{ak}{n}}\ra g'(v).$$
   Consequently,
   $$\frac{c_3(u-a, v)-c_3(u,v)}{a}+g'(v)=0.$$
 Hence $c_3(u-a,v)-c_3(u,v)=-a g'(v)$ for all $u\in \R$ and all
 $a>0$.  It follows that
  $c_3(u,v)=c_3(0, v)+g'(v) u$ for all $u$.

 \end{proof}

\b{Le}\label{l5.3}
 Suppose $g$ is differentiable at $v_1$ and $v_2$.  Then
 $g'(v_1)=g'(v_2)$.

 \end{Le}

 \b{proof}  By Lemma \ref{l5.2}, we have
  $c_3(u, v_1)=c_3(0, v_1)+u g'(v_1)$ and
   $$c_3(u+[v_2-v_1]\ln|v_2-v_1|, v_2)=c_3(0,
   v_2)+(u+[v_2-v_1]\ln|v_2-v_1|)g'(v_2)$$
       for all $u$.
Now   Lemma \ref{l5.1}  applied  to $u, v_1$, $v_2$ implies  that
$\big\vert u(g'(v_2)-g'(v_1))\big\vert\le C$
  holds for all $u$, where $C$ is a quantity independent  of $u$.
 It follows that $g'(v_2)-g'(v_1)=0$.

 \end{proof}

\noindent
    {\bf{Completing the proof of Theorem \ref{Jordan}}}.
        Lemma \ref{l5.3} implies that $g$ is an affine
function and hence there are constants $a, b$ such that $g(v)=av+b$.
It now follows from Lemma \ref{l5.2} that  for any $v$ we have
$c_3(u,v)=c_3(0,v)+au $. To finish the proof of Theorem
\ref{Jordan}, it remains to show that $c_3(0, v)$ is Lipschitz in
$v$.  This follows immediately from Lemma \ref{l5.1} after plugging
in the formulas for $g$ and $c_3$.

 Now the proof  of Theorem \ref{Jordan} is complete.

 \qed

\section{A Liouville type theorem }\label{liouv}

In this section we  prove  a Liouville type theorem for $G_A$ in the
case when $A$ is a Jordan block: every conformal map of the ideal
boundary of $G_A$ extends to an isometry of $G_A$.

Let  $X$ and $Y$  be  quasimetric spaces with finite Hausdorff
dimension. Denote by $H_X$ and $H_Y$  their Hausdorff dimensions and
by $\mathcal{H}_X$  and $\mathcal{H}_Y$  their Hausdorff measures.
   %(see
  %\cite{F}   for definitions).
       We say  a quasisymmetric map $f: X\ra Y$ is
  conformal if:\newline
  (1) $L_f(x)=l_f(x)\in (0, \i)$  for $\mathcal{H}_X$-almost every
  $x\in X$;\newline
  (2) $L_{f^{-1}}(y)=l_{f^{-1}}(y)\in (0, \i)$ for
  $\mathcal{H}_Y$-almost every $y\in Y$.

We now  describe some isometries of $G_A$. For any $g=(x,t)\in
G_A=\R^n\rtimes \R$, the Lie group left translation $L_g$ is an
isometry.  If $g=(x,0)$,  then  the boundary map
  $\p L_g: \R^n\ra \R^n$ of $L_g$ is
 translation by $x$.  If $g=(0,t)$, then the boundary map of $L_g$
 is the similarity $e^{tA}$.  Let $\tau':  G_A\ra G_A$ be given by
 $\tau'(x,t)=(-x,t)$. Then $\tau'$ is an isometry, and its boundary
 map
 is  $\tau: \R^n\ra \R^n$, $\tau(x)=-x$.

\b{Th}\label{liu}
  Let $J_n$ be the $n\times n$ ($n\ge 2$) Jordan matrix with
    eigenvalue $1$. Then every conformal map $F: (\R^n, D_{J_n})\ra (\R^n,
D_{J_n})$ is
  the boundary map of an isometry  $G_{J_n}\ra G_{J_n}$.

\end{Th}

\b{proof}  Let $F: (\R^n, D_{J_n})\ra (\R^n, D_{J_n})$  be  a
quasisymmetric map. After composing with the boundary maps of
isometries described above, we may assume $F$ has the following form
   $$F(x)=(I+a_1N+\cdots +a_{n-2}N^{n-2})x+(C(x_n), 0,
   \cdots, 0)^T,$$
    where $C:\R\ra \R$ is Lipschitz.
  %By post-composing with an Euclidean translation, we may further
  %assume $C(0)=0$.
 We will prove the following statement by inducting on $n$:\newline
  $$\text{\e{If $F$ as above is conformal, then $a_1=\cdots=a_{n-2}=0$ and
   $C$  is constant.}}$$

  The basic step
$n=2$ is Theorem 6.3  in \cite{X}.
 Now we assume $n\ge 3$ and that the statement holds for Jordan
 matrices with sizes $\le n-1$.
     Notice that  $F$ maps every  horizontal hyperplane $H(x_n):=\R^{n-1}\times \{x_n\}$
    to itself.
       By Lemma \ref{l11}  the restriction of $D_{J_n}$  on $H(x_n)$  agrees with the metric
    $D_{J_{n-1}}$.
  It now follows from Fubini's theorem that for a.e. $x_n\in \R$,
      the restricted map
    $$F|_{H(x_n)}: (H(x_n),
    D_{J_{n-1}})\ra  (H(x_n),
    D_{J_{n-1}})$$
     is also conformal.
  %We may assume $x_n=0$ after composing $F$ with Euclidean
  %translations.
     Now the induction hypothesis
    applied to $F|_{H(x_n)}$ implies that $a_i=0$ for $1\le i\le n-2$.
     It remains to show $C$  is constant.

     Suppose $C$ is not constant.   Then  there is some $u\in \R$ such that  $C'(u)\not=0$
       and   $L_F(p)=l_F(p)$ for some $p\in H(u)$.  After pre-composing and post-composing with Euclidean
     translations, we  may assume  $u=0$,  $C(0)=0$   and $p$ is the origin $o$.
       %$C'(0)\not=0$ and $C(0)=0$.
Notice that the restriction of $F$ to the $x_1$-axis is the
identity, so
 $L_F(o)=l_F(o)=1$.  Now for any $x_n>0$, choose $x_1, \cdots,
 x_{n-1}$ such that $x=(x_1, \cdots x_n)^T$ satisfies
   $e^{-tN}x=(0, \cdots, 0, x_n)^T$, where $t=\ln x_n$.
    It follows that $D_{J_n}(o, x)=e^t=x_n$.
 Suppose $D_{J_n}(F(o), F(x))=e^s$. Then
  $e^s=|e^{-sN} F(x)|$.  We calculate as before that
   $$e^{-sN}F(x)=\left(C(x_n)+\frac{(t-s)^{n-1}}{(n-1)!}x_n, \;
   \frac{(t-s)^{n-2}}{(n-2)!}x_n,\;
   \cdots, \; (t-s)x_n,\; x_n\right)^T.$$
       Since   $L_F(o)=l_F(o)=1$,
     we must have
     $\frac{e^s}{e^t}=\frac{D_{J_n}(F(x), F(o))}{D_{J_n}(x, o)}\ra 1$  as $x_n\ra 0$  and
     hence $t-s\ra 0$.
Now
$$e^s=\big\vert e^{-sN} F(x)\big\vert=x_n \bigg\vert\left(\frac {C(x_n)}{x_n}
+\frac{(t-s)^{n-1}}{(n-1)!},\;\frac{(t-s)^{n-2}}{(n-2)!},\;\cdots,\;
(t-s),\; 1\right)^T\bigg\vert.$$
 Since $x_n=e^t$, we have
 $$e^{s-t}=\bigg\vert\left(\frac {C(x_n)}{x_n}
+\frac{(t-s)^{n-1}}{(n-1)!},\;\frac{(t-s)^{n-2}}{(n-2)!},\;\cdots,\;
(t-s),\; 1\right)^T\bigg\vert.$$
     Now  as $x_n\ra 0$,   the right hand side converges to
     $$\big\vert\left(C'(0), 0, \cdots, 0,
1\right)^T\big\vert=\sqrt{1+(C'(0))^2},$$
  which is $>1$  since $C'(0)\not=0$.  However, the left hand side converges to $1$.  The
 contradiction shows $C$  must be  a constant function.

\end{proof}

 \addcontentsline{toc}{subsection}{References}

\noindent Address:

\noindent Xiangdong Xie: Dept. of Mathematical Sciences, Georgia
Southern University, Statesboro, GA 30460, U.S.A.\hskip .4cm E-mail:
xxie@georgiasouthern.edu

\end{document}